\documentclass[symmetry,article,submit,oneauthor,
pdftex,12pt,a4paper]{mdpi}

\usepackage{amsfonts,amssymb}
\pdfoutput=1
\lastpage{x}
\doinum{10.3390/------}
\pubvolume{xx}
\pubyear{2012}
\history{Received: xx / Accepted: xx / Published: xx}

\newcommand{\lt}{\triangleleft}

\newcommand{\Co}{\mathbb{C}}

\newcommand{\R}{\mathbb{R}}

\newtheorem{theorem}{Theorem}[section]

\usepackage{graphicx}

\Title{Classical Knot Theory}

\Author{J. Scott Carter $^{1,\star}$}

\address{%
$^{1}$ University of South Alabama, Department of Mathematics and Statistics, ILB 325, Mobile, AL, USA}

\corres{E-mail carter@jaguar1.usouthal.edu, Phone:1-251-460-6264}

\abstract{ This paper is a very brief introduction to knot theory. It describes knot coloring by quandles, the fundamental group of a knot complement, and handle-decompositions of knot complements.}


\begin{document}

\section{Introduction}
Over the past several months, Slavik Jablan asked me to contribute to this volume. I hesitated since I have very little to say that is new, that which is new does not fit into this venue, and I have been procrastinating rather than writing the newer results. 
 The current article contains one of those results: the handle decomposition that yields the Alexander-Briggs \cite{AB} presentation of the knot group. Masahico Saito, Dan Silver, Susan Williams, and I will be writing more about this presentation as it relates to virtual knots and knot colorings shortly after this article is complete. The article that I present to you now could be written by any other author in the current volume. Many could do a better job than I, for I am focusing upon some well-known properties of classical knots and the space that surrounds them. 

Specifically, the article focuses upon the so-called handle decomposition of the knot complement. These handle-theoretic techniques are central to geometric topology. My hope is that the details presented here will aid an uninitiated reader in mastering them. 

This paper was developed as an introductory series of lectures that I gave to the beginning graduate students at Kyongpook National University, Daegu, Korea during a sabbatical semester from the University of South Alabama. My visit in Korea is being funded by a grant from the Brain-Pool Trust. I would like to thank my host, Professor Yongju Bae for his hospitality. I also would like to thank my long-time collaborator Masahico Saito for our discussions about these matters and related things. Finally, many thanks to the referees for their helpful and kind comments.

Here is an outline. The second section follows this introduction. In it,   I indicate that the three knots that are illustrated in Fig.~\ref{untreate} are distinct. The discussion there focuses upon Reidemeister moves, quandle colorings, and the quandle structure of a group under the conjugation operation. The third section gives an outline of the definition of the fundamental group of a knot complement. I have included some of the diagrams that are used to demonstrate the homotopy equivalences necessary for the well-definedness of the invariant. Many texts ({\it e.g.} \cite{CF,Hatcher,Chuck}) present this material in detail. The fourth section sketches handle theory. It describes handle decompositions of a few surfaces (sphere, torus, projective plane, Klein bottle) illustrates handle sliding in two dimensions and discusses turning the handle decomposition of the torus upside-down. In the fifth section, the discussion turns to decompositions of the trefoil and the figure-eight knot complements. The knot diagrams  are annotated to manipulate these decompositions. In the last section, I describe an alternative handle decomposition that can be used to present the fundamental group of a knot. I call this presentation the Alexander-Briggs presentation since a careful reading of the last section their seminal paper  \cite{AB} indicates that these authors understood this handle decomposition even if handles had not yet been defined during that era. 

This paper does not cover any post-Jones ({\it e.g.} \cite{Jones,homfly}) invariants and it does not cover the Alexander polynomial even though the latter can easily be gleaned from some of the discussions herein. For a wonderful historical survey see \cite{Przytycki}.  I am writing with a deadline in mind. My goal is a self-contained opus that contains the ideas and visual imagery that occupies my current state of consciousness. 

Let us begin.
%



%

\section{The unknot, the trefoil, and the figure-8 knot are distinct}

\begin{figure}[htb]
\begin{center}
\includegraphics[width=3in]{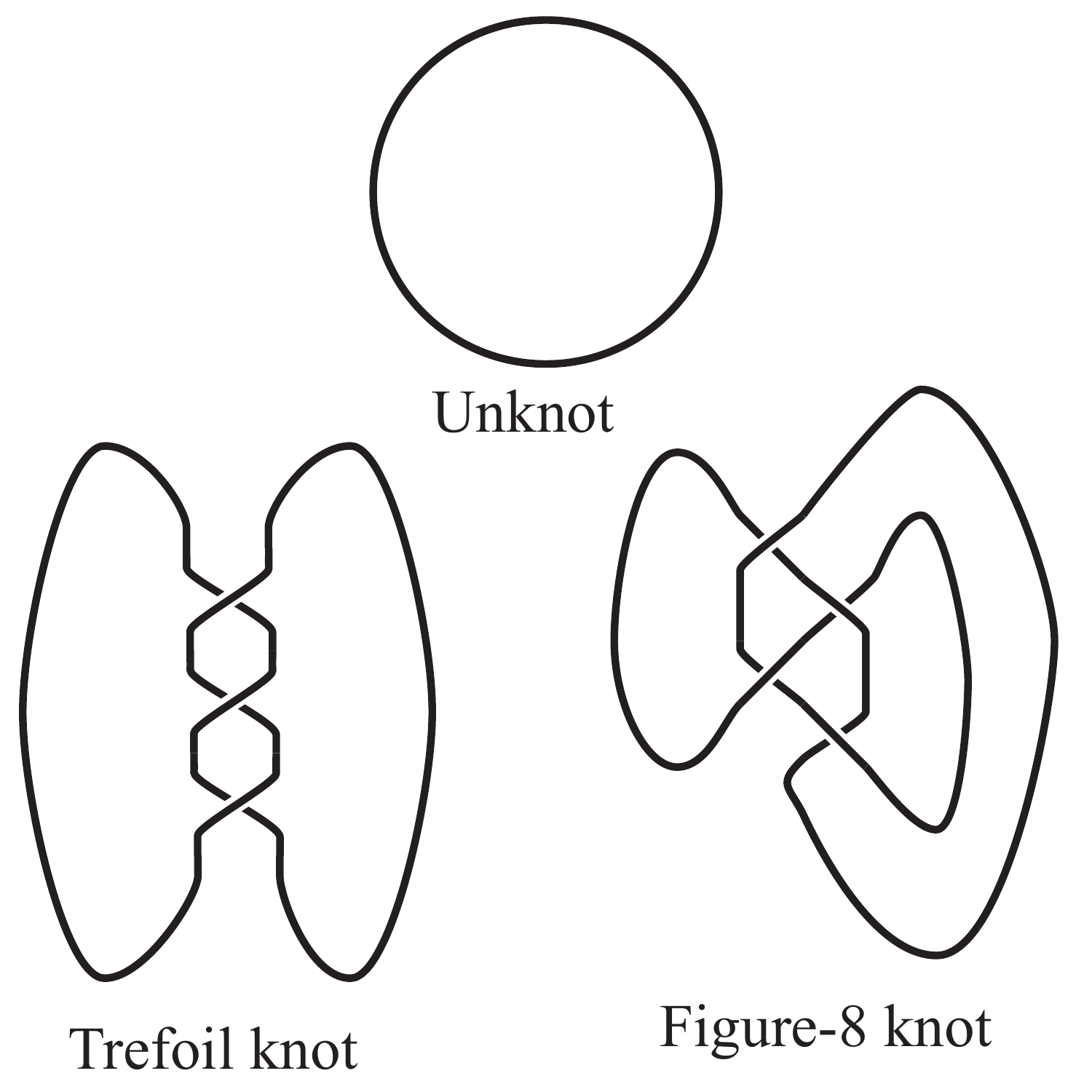}
\end{center}
\caption{The unknot, the trefoil knot, and the figure-8 knot}
\label{untreate}
\end{figure}

The goal of this section will be to demonstrate that the three knots that are depicted in Fig.~\ref{untreate} are distinct. 
The embedded curves represented are called {\it  the  unknot, the trefoil knot, and the figure-8 knot} as indicated in the figure. First, here are some basic notions.

A {\it classical knot} is a (smooth or piecewise-linear locally-flat) embedding of a circle $S^1 =\{ z\in \Co: |z|=1 \}$ into $3$-dimensional space. Two such knots are said to be equivalent if one can be continuously deformed into the other without breaking or cutting. More precisely, $f_1:S^1 \hookrightarrow \R^3$ and $f_2:S^1 \hookrightarrow \R^3$ are {\it equivalent} if and only if there is an orientation preserving homeomorphism of pairs $(\R^3, f_1(S^1)) \rightarrow (\R^3, f_2(S^1))$. We often replace $\R^3$ with $S^3 = \{ (x,y,z,w): x^2 +y^2 + z^2 + w^2 =1 \}$ which is the one-point compactification of $3$-space.

 \begin{figure}[htb]
\begin{center}
\includegraphics[width=2in]{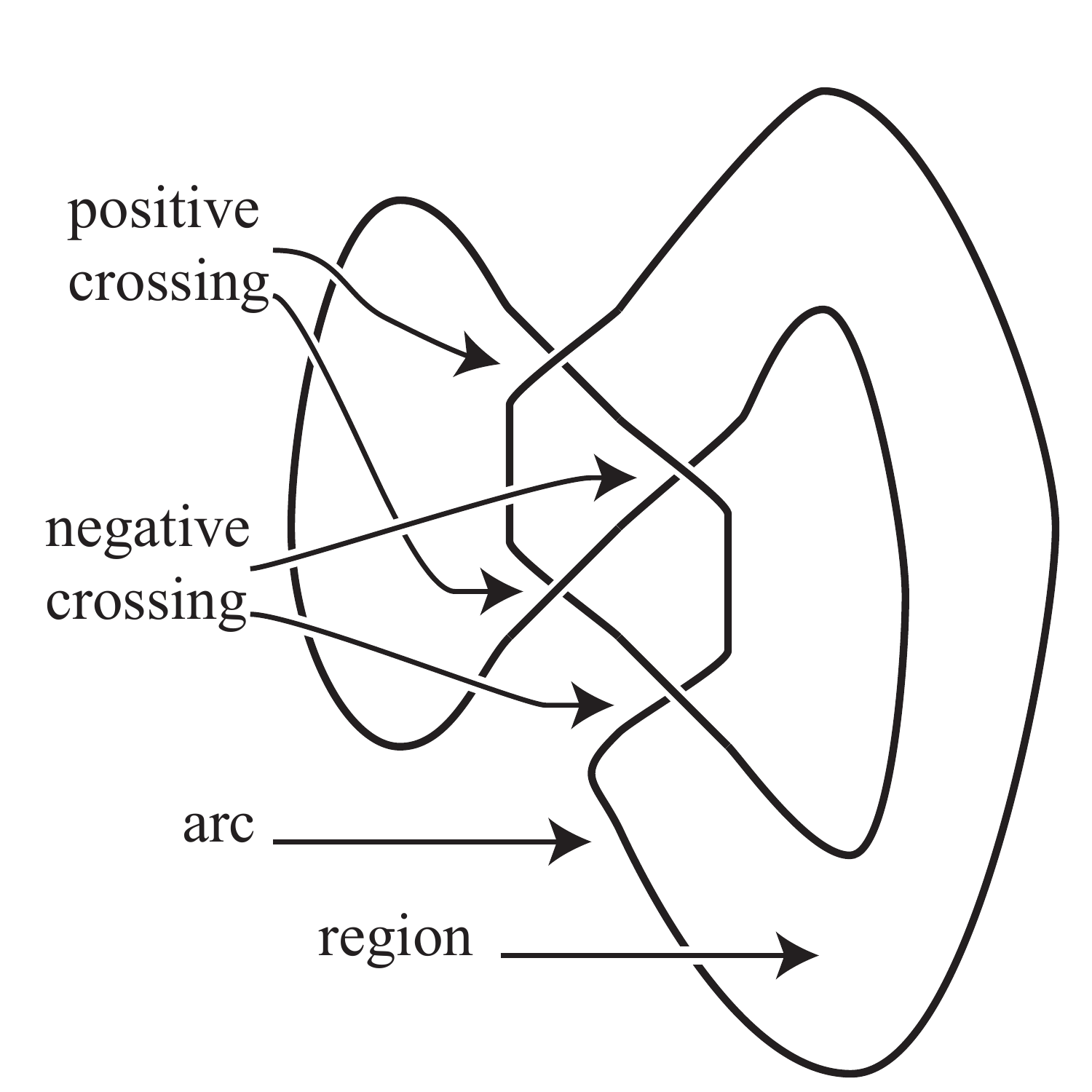}
\end{center}
\caption{A diagram of the figure 8 knot}
\label{fig8}
\end{figure}

\begin{figure}[htb]
\begin{center}
\includegraphics[width=2in]{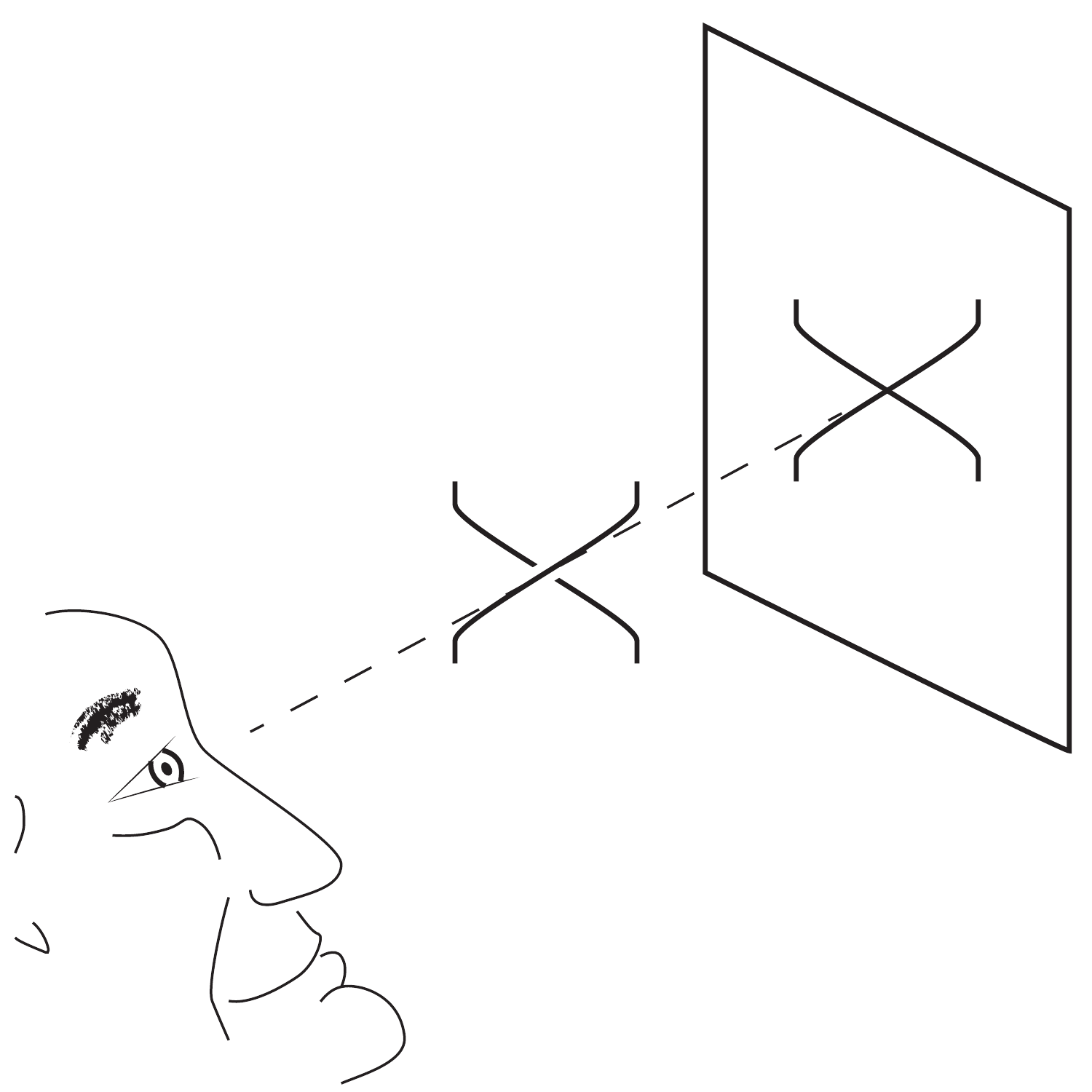}
\end{center}
\caption{The projection of a crossing and point of view}
\label{eyeplane}
\end{figure}

A classical knot is depicted via {\it a diagram} as indicated in Fig.~\ref{fig8}. In such a diagram a generic projection of the knot to plane is chosen; when two arcs project to transversely intersecting arcs in the the plane, crossing information is depicted by breaking the arc that is closer to the plane of projection than an observer is. Some additional standard terminology is indicated in that figure. Figure~\ref{eyeplane} indicates crossing conventions. From a diagram, a local picture of a crossing can be reproduced via choosing a coordinate system at the intersection points and lifting the under-arc (which lies along the $y$-axis) to the segment $\{(0,y,1): |y|\le 1\}$, and lifting the over-arc to the segment 
$\{ (x,0,2): |x| \le 1\}$. Observe that if the orientations of these segments coincide with the orientations of the axes, then the crossing so depicted is positive.

Knot diagrams are studied via an equivalence relation on diagrams. The relation is generated by the local moves that are indicated in Fig.~\ref{reid3} and which are called {\it Reidemeister moves}. These moves reflect the motion of a knot in space in the sense that if a knot is generically moved through space, then the projection of the motion is encapsulated by a sequence of applications of moves to the projection. In order to experience the moves directly, I suggest making a knot out of a rigid material such as thick wire and watch the knot as you move your point of view. It is not difficult to experience each of the Reidemeister moves in one's visual sphere.

\begin{figure}[htb]
\begin{center}
\includegraphics[width=2.5in]{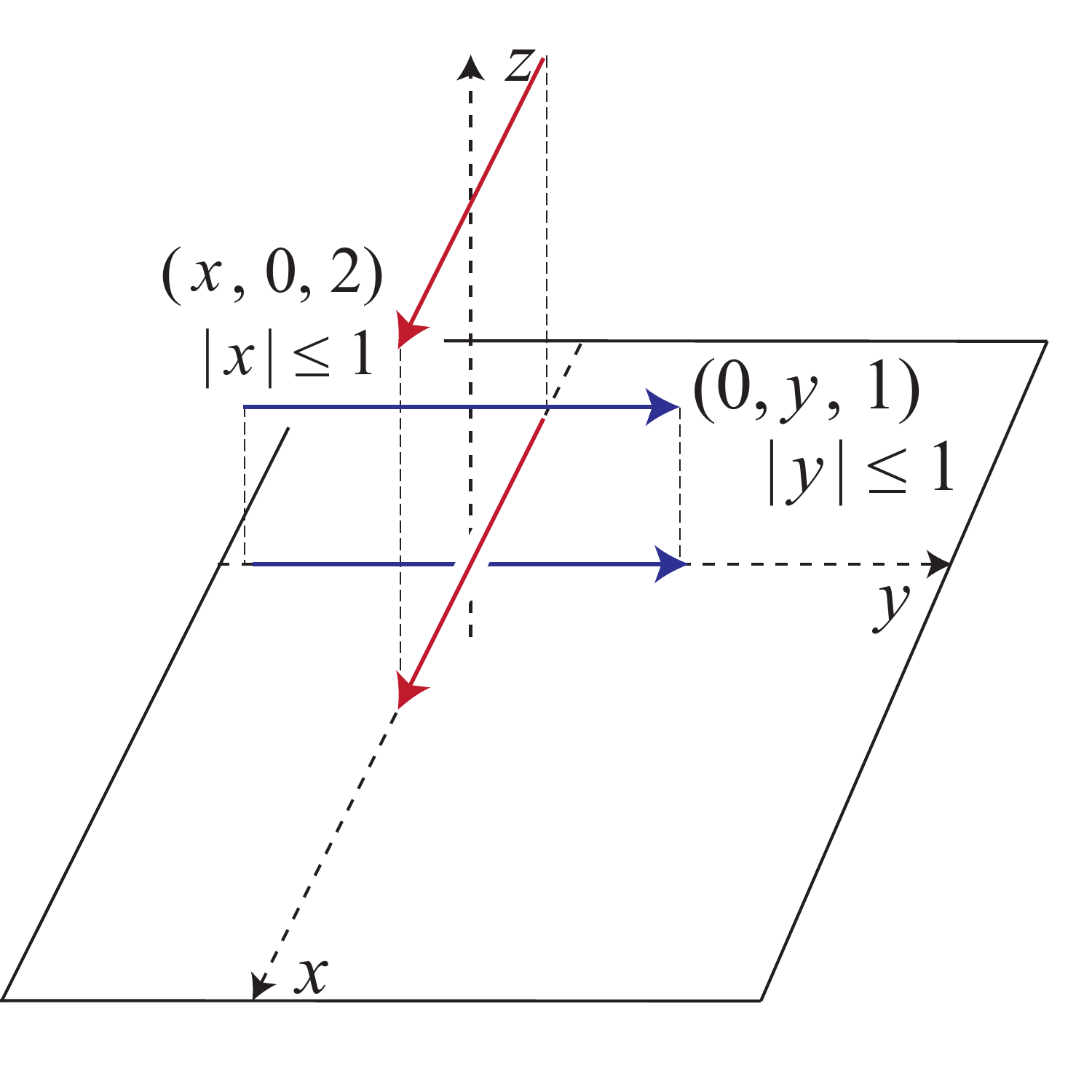}
\end{center}
\caption{Lifting a crossing to standard position in space}
\label{liftingcrossing}
\end{figure}

With these preliminary matters established, we demonstrate that the trefoil is not the unknot. To do so, we color the arcs of the trefoil using three distinct colors as indicated in Fig.~\ref{tre3}. In general a knot is said to be {\it 3-colorable} if each of the arcs in a representative diagram can be assigned a color such that either all three colors are coincident at a crossing or only one color is incident.
This relationship is indicated in Fig.\ref{tri}.

\begin{figure}[htb]
\begin{center}
\includegraphics[width=4in]{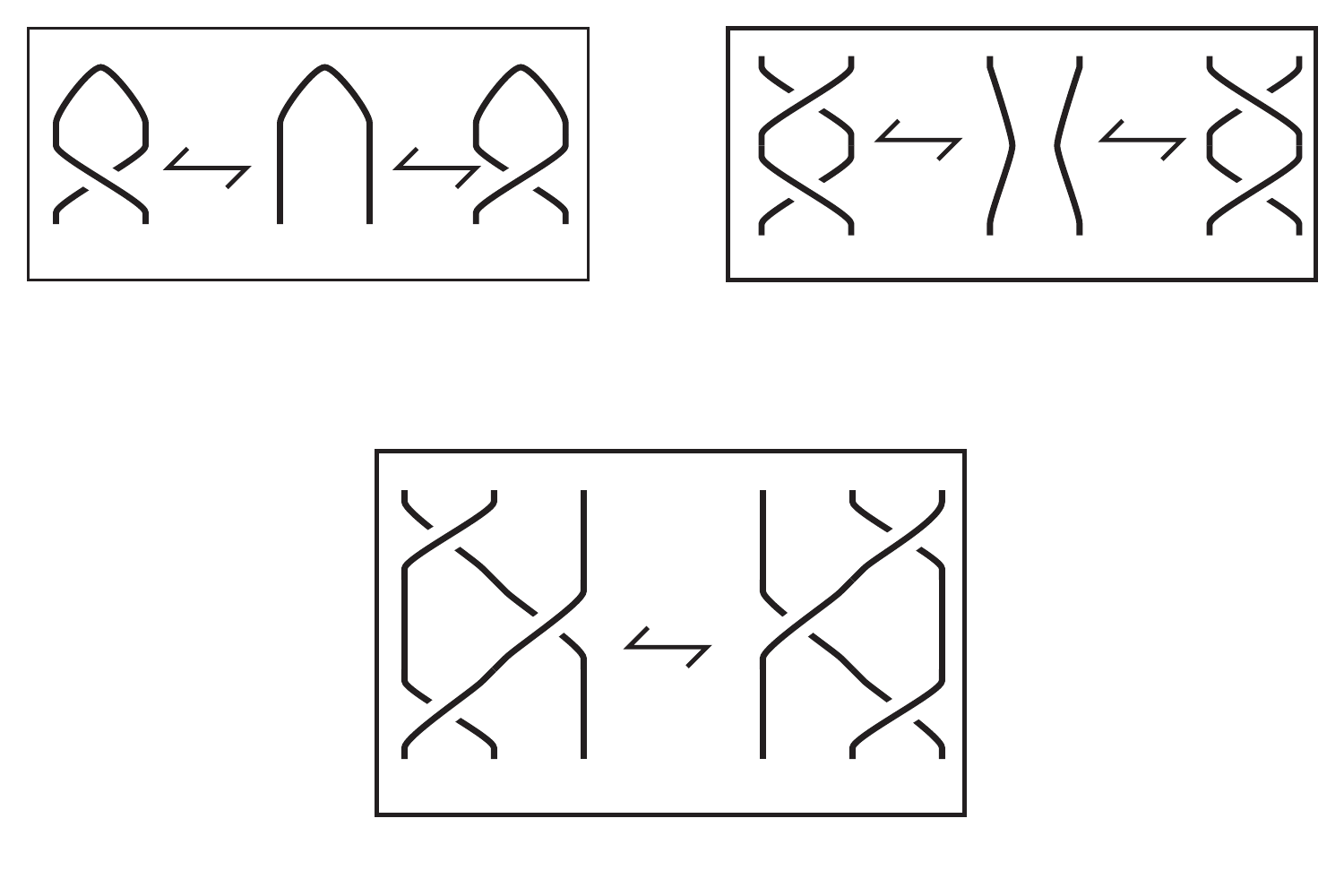}
\end{center}
\caption{The Reidemeister moves}
\label{reid3}
\end{figure}

\begin{figure}[htb]
\begin{center}
\includegraphics[width=2in]{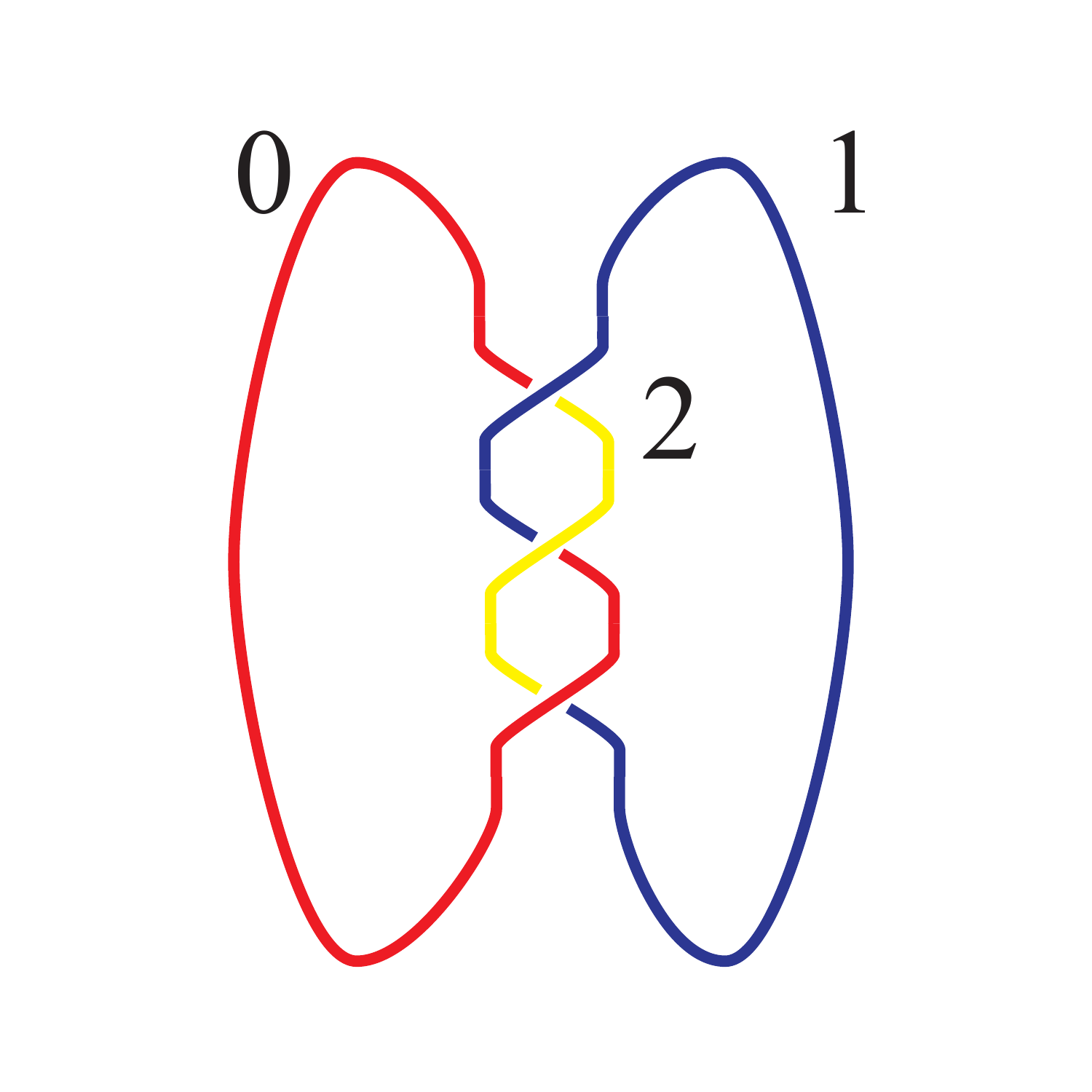}
\end{center}
\caption{The trefoil knot is 3-colorable}
\label{tre3}
\end{figure}

\begin{figure}[htb]
\begin{center}
\includegraphics[width=1.5in]{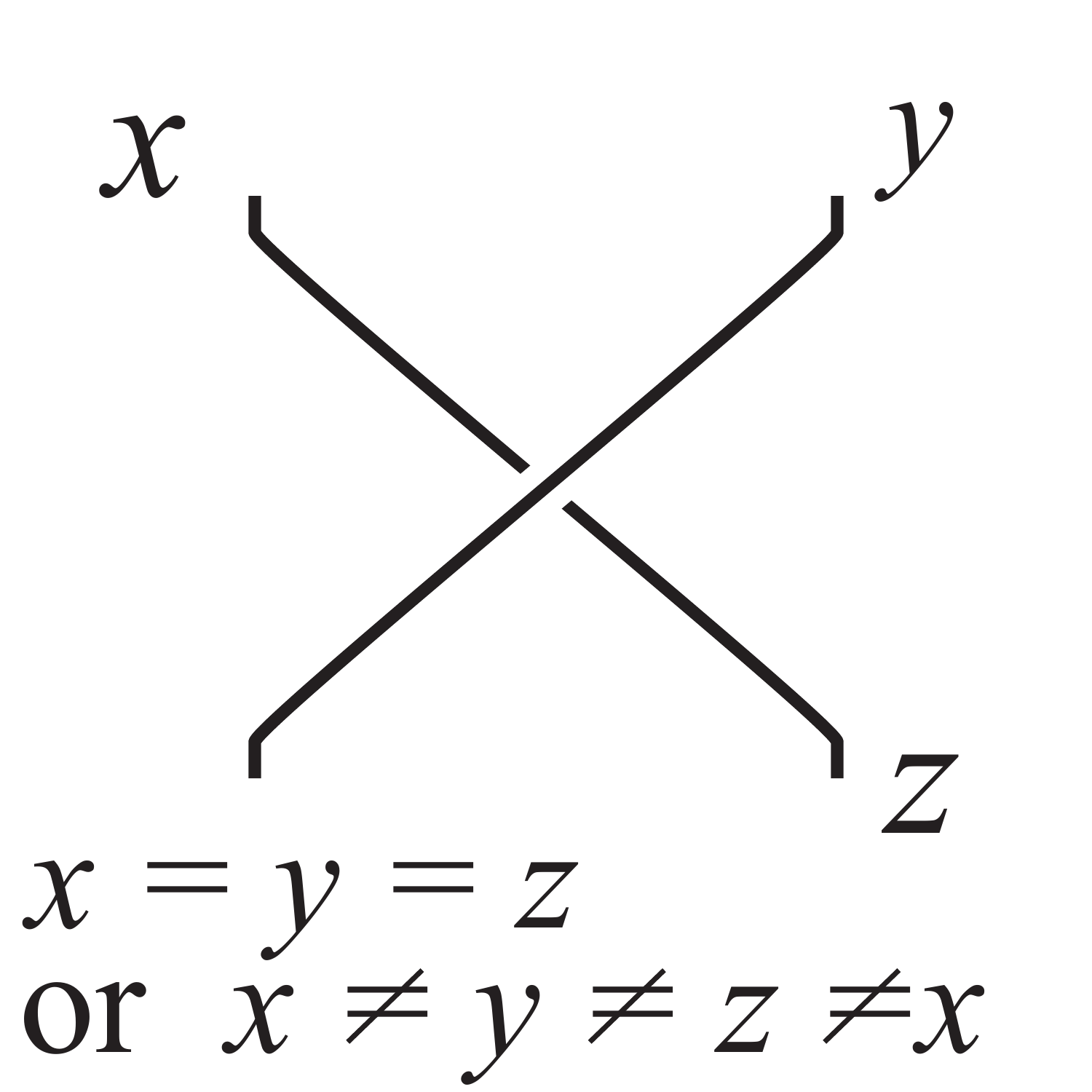}
\end{center}
\caption{The 3-colorabilty condition}
\label{tri}
\end{figure}

Since $3$-colorability is defined in terms of a representative diagram, we must show that it is invariant under the Reidemeister moves. Observe that for a type I move (top of Fig.~\ref{reid3}) all three colors must coincide. For a type II move (middle of the figure), either the two disentangled arcs have different colors or they have the same color. In the latter case, the short arc that is introduced has the same color as these two arcs do. In the former case, the short arc has the third color.

For the type III move there are five cases to consider:
\begin{enumerate}
\item Case $(a,a,a)$: In this case, the  three arcs from left to right on the left-hand-side of the move all have the same color at the top of the diagram. All the arcs on either side of the move are the same color.

\item Case $(a,b,c)$: the three arcs from left to right on the left-hand-side of the move all have different colors at the top of the diagram. The top of the arcs on the right hand side are similarly colored. At the bottom of the diagram on either side the colors are  $(c,a,c)$ from left to right.

\item Case $(a,a,b)$: the first two arcs at the top of either the left or right side of the move have the same colors  and the last arc is colored differently.  At the bottom of either side the colors are $(b,c,c)$ from left to right.

\item Case $(a,b,b)$: The two arcs on the right have the same color,  and the first arc is colored differently. 

\item Case $(c,a,c)$: This case is obtained from case $(a,b,c)$ by turning it upside-down.  

\end{enumerate}

\begin{figure}[htb]
\begin{center}
\includegraphics[width=3in]{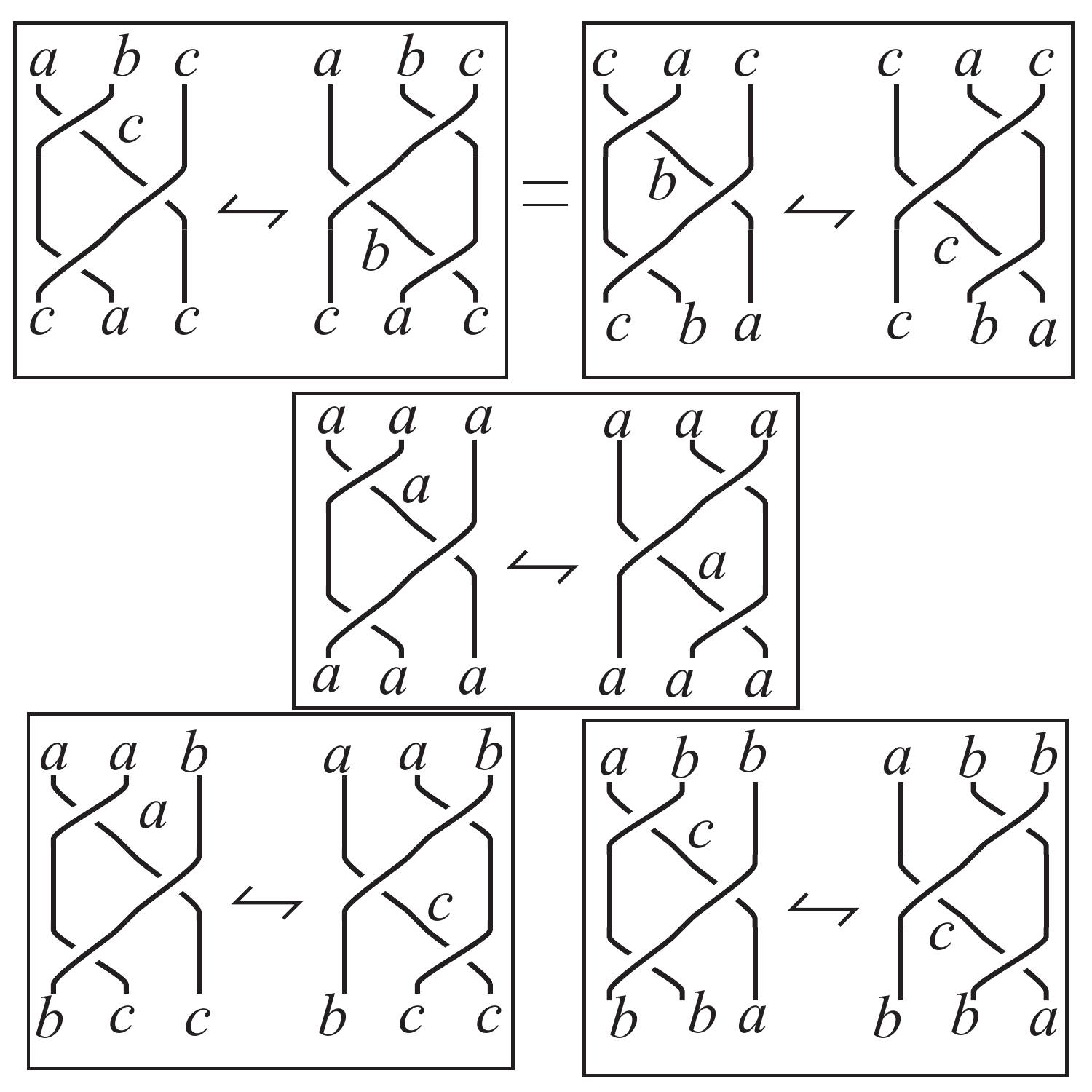}
\end{center}
\caption{The type III move and colorability}
\label{colorIII}
\end{figure}

Figure~\ref{colorIII} illustrates.

From these illustrations, we can see that if a diagram is $3$-colorable, then any diagram of the knot is $3$-colorable in a way that uses three distinct colors. Moreover, any coloring of the unknot is monochromatic. Thus the unknot is distinct from the trefoil.  Any attempt to $3$-color the figure-$8$ knot with three distinct colors will fail. We leave that as an exercise for the reader. 

To demonstrate that the figure-$8$ knot is distinct from the other examples, we formalize the idea of colorability be demonstrating that crossings in a knot diagram can be used to axiomatize an algebraic system: a set with a binary operation that is assumed to satisfy three axioms which correspond to the Reidemeister moves. Specifically, if an over-crossing arc is oriented, then  a homunculus standing on the over-arc holding out its{\footnote{all homunculi in this paper are gender neutral}} left hand points towards an under-arc. The under-arc on the right is labeled $a$, the over-arc is labeled $b$, and the target arc is labled with a product $a \lt b$. We read  $a\lt b$ as ``$a$ acted upon by $b$."


\begin{figure}[htb]
\begin{center}
\includegraphics[width=2in]{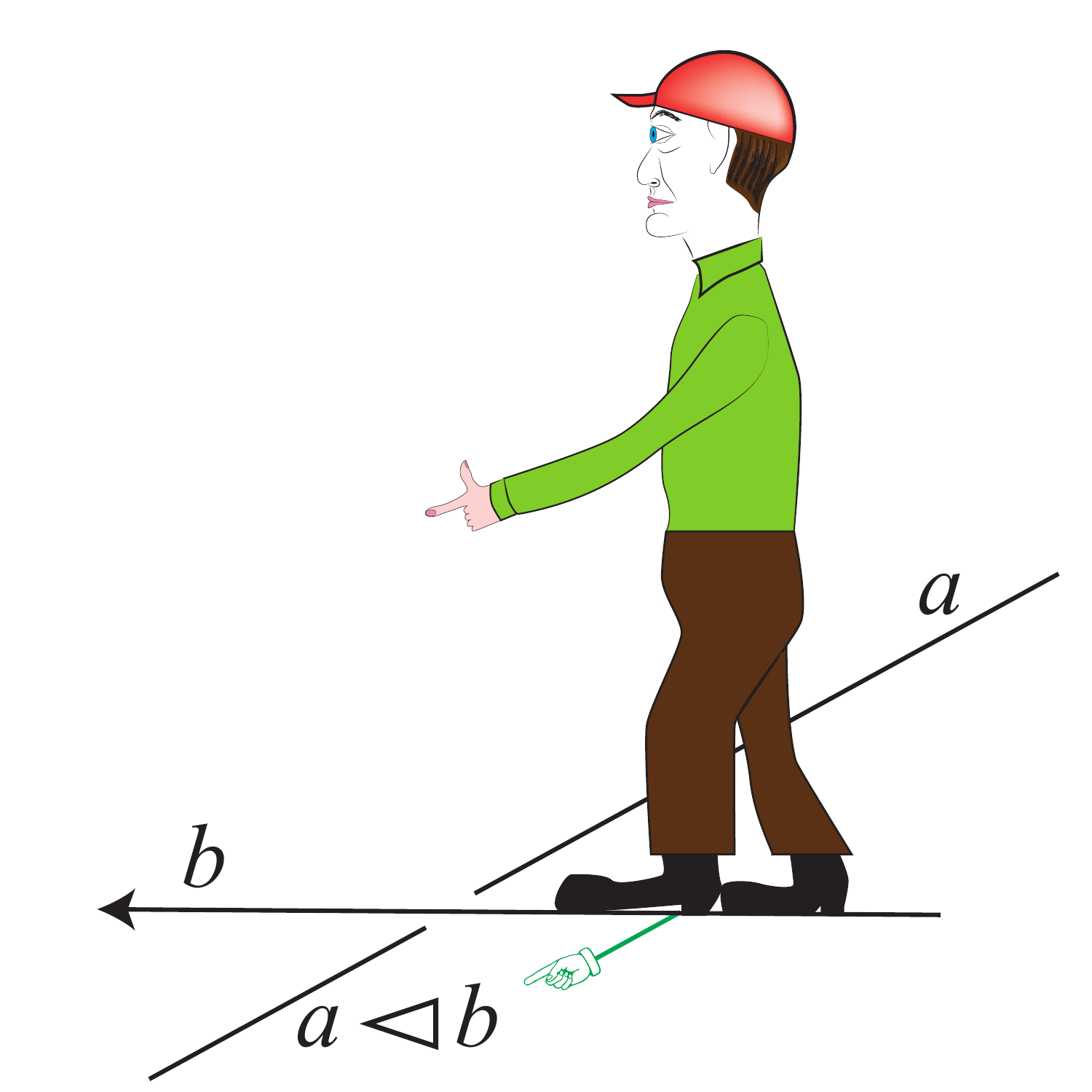}
\end{center}
\caption{The under-arc towards which the left hand points received the product}
\label{pointingleft}
\end{figure}

\begin{figure}[htb]
\begin{center}
\includegraphics[width=2.5in]{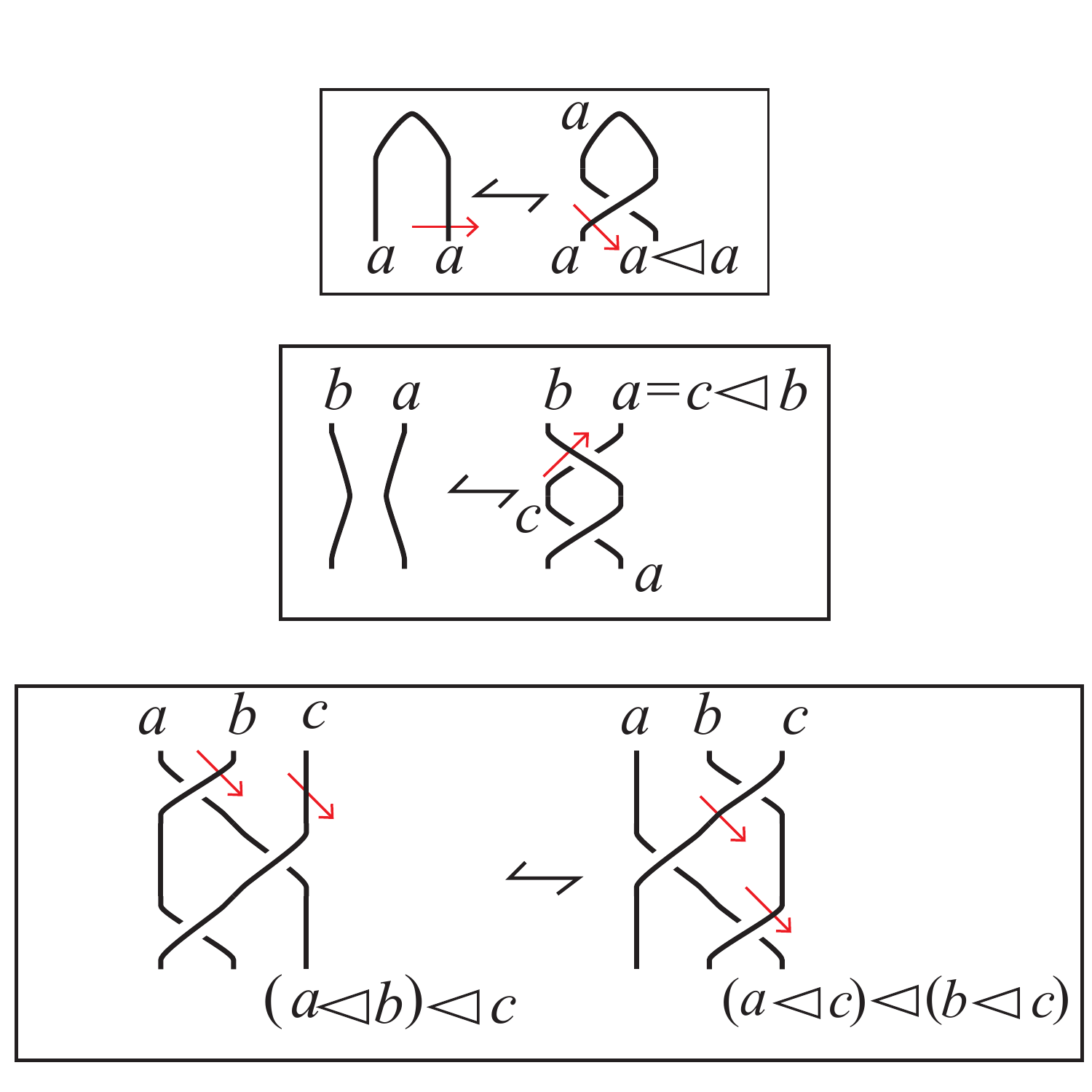}
\end{center}
\caption{The axioms of a quandle correspond to the Reidemeister moves}
\label{reidQ}
\end{figure}

To continue, we define a set $X$ that has a binary operation $\lt: X\times X \rightarrow X$ (written in in-fix notation) a {\it quandle} if the following three axioms hold:
\begin{enumerate}\item[I] for any $a\in X$, we have $a\lt a =a$;
\item[II] for each $a,b \in X$, there is a unique $c\in X$ such that $c\lt b = a;$
\item[III] for each $a,b,c \in X$, we have $(a\lt b)\lt c = (a\lt c)\lt (b \lt c).$
\end{enumerate}
As indicated in Fig.~\ref{reidQ}, the axioms correspond to the Reidemeister moves in the sense that if a knot diagram is colored by elements of a quandle $X$ (in such a way that the under-arc towards which the humunculus's left hand points receives a product), then any diagram related by a sequence of Reidemeister moves will also be colored. 

The first example of a quandle that we have in mind is the set $X=\{0,1,2\}$ with $a\lt b = 2b -a  \  ({\mbox{\rm mod}} \ {3})$. This product yields the $3$-colorability conditions of the previous example. The quandle here is called $R_3$ --- the dihedral quandle of order $3$.

If $X$ is a subset of a group that is invariant under conjugation by its elements, then $X$ is a quandle under the operation $x\lt y = y^{-1} x y$. For example, let $G= S_4 =S_4\{0,1,2,3\}$ denote the group of permutations of the set $\{0,1,2,3\}$. Let $X=\{(1,2,3),(0,3,2),(0,1,3),(0,2,1)\}$ denote the subset of ``oriented" $3$-cycles. Label each such  $3$-cycle by the element that it fixes: $0 \leftrightarrow (1,2,3)$, $1\leftrightarrow (0,3,2),$ $2\leftrightarrow (0,1,3)$, and $3\leftrightarrow (0,2,1)$. Then the quandle table for $X$ is 
as follows: 
\begin{center} 
\begin{tabular}{||c||c|c|c|c||}\hline \hline row $
\lt$ col & 0 & 1 & 2 & 3 \\ \hline \hline
0&0&3&1&2 \\ \hline
1 & 2 & 1 & 3& 0 \\ \hline
2 & 3 & 0 & 2 &1 \\ \hline
3 & 1 & 2 & 0 & 3 \\ \hline \hline \end{tabular}\end{center}

\begin{figure}[htb]
\begin{center}
\includegraphics[width=2.5in]{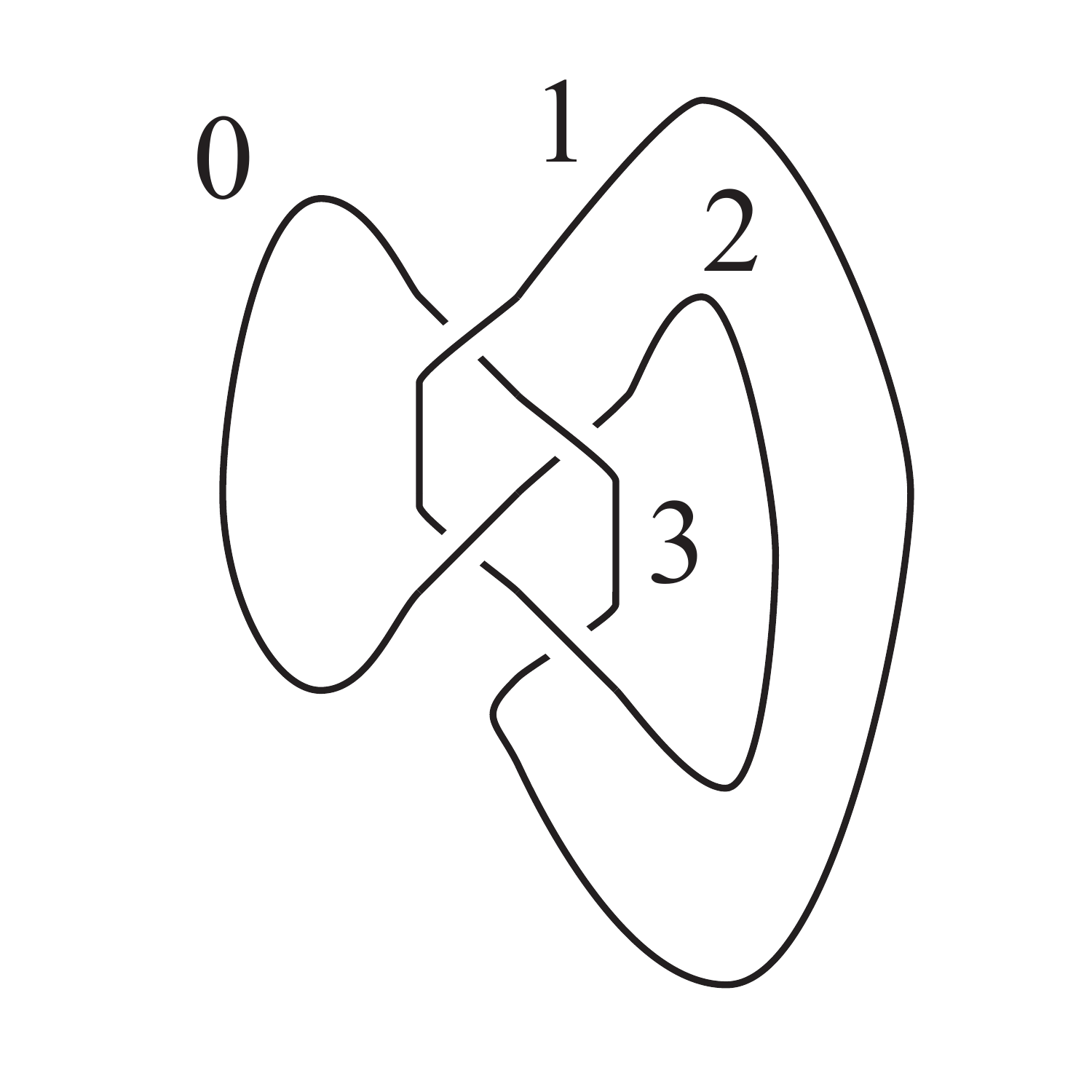}
\end{center}
\caption{The figure-8 knot can be colored by $QS_4$}
\label{eghtbyfour}
\end{figure}

This quandle is called ${\mbox{\rm QS}}_4$ --- {\it the tetrahedral quandle}.
Figure~\ref{eghtbyfour} indicates that the figure-8 knot can be colored by ${\mbox{\rm QS}}_4$. Thus the figure-8 knot is distinct from the unknot. Even though the trefoil knot can be colored by ${\mbox{\rm QS}}_4$, it is distinct from the figure-8 knot since the latter cannot be colored by $R_3$. 

We state the classical result:

\begin{theorem} The unknot, the trefoil knot, and the figure-8 are distinct. \end{theorem}

In fact, finite quandles are very good at distinguishing knots. Each classical knot has an associated fundamental quandle that  is somewhat stronger than its fundamental group. The development of the quandle as a knot invariant occurred independently in papers  by Matveev~\cite{Matveev} and Joyce~\cite{Joyce} about twenty years before the turn of the century. 

The dihedral quandle and the tetrahedral quandle are both examples of quandles that are defined geometrically. The dihedral quandle $R_n$ consists of the reflections of an $n$-gon composed under conjugation. The tetrahedral quandle is defined similarly as rotations of the tetrahedron through a vertex. In general, a (connected) quandle exhibits aspects of symmetry since it can be defined in terms of a binary operation on a set of cosets of a group of automorphisms of the quandle. This description appears elsewhere \cite{JSCQ,Joyce,Matveev} and would take us far away from the current purposes of the paper.

\section{The fundamental group}

Let $K:S^1\hookrightarrow S^3$ denote a smooth or PL-locally-flat embedding of a circle into the {\it $3$-dimensional sphere} $S^3=\{ (x,y,z,w) \in \R^4: x^2+y^2+z^2+w^2 =1\}$. For brevity, we write the image of the embedding as $K$, and we speak of the {\it knot $K$.} The smooth or PL-locally-flat condition suffices to provide a tubular neighborhood of the knot. This is a smooth embedding of a solid torus $N: S^1 \times D^2 \hookrightarrow S^3$ that is a {\it tubular neighborhood} of the knot. That is, letting $D^2 =\{ (x,y)\in \R^2: x^2 +y^2 \le 1 \}$ denote the {\it $2$-dimensional disk}, then the knot is embedded as the core $\{0\} \times S^1$ of this solid torus. 
The {\it knot exterior} is the space $E=E(K)=S^3 \setminus {\mbox{\rm int}}(N).$ It has a torus ($S^1\times S^1$) as its boundary. We fix a base point of the exterior, and call the point $*$. \begin{figure}[htb]
\begin{center}
\includegraphics[width=5in]{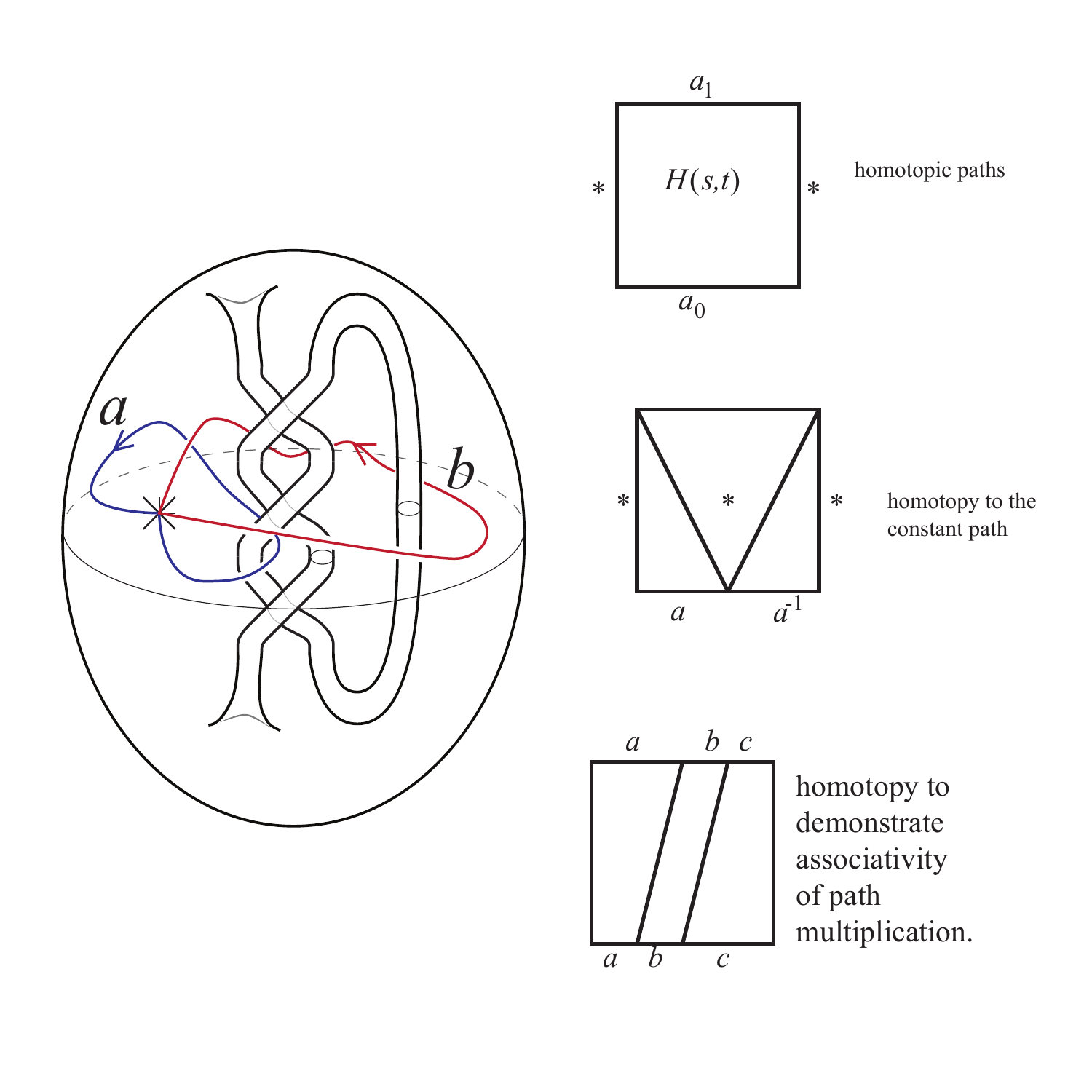}
\end{center}
\caption{Ideas that are used to demonstrate the fundamental group and necessary homotopies}
\label{fundy}
\end{figure}

The {\it fundamental group of the knot exterior $E$  based at the point $*$} is defined as the set of homotopy classes of maps of pairs $\gamma: ([0,1],\{0,1\})\rightarrow (E,*)$ where the homotopies are required to fix the boundary points
at the base point, and the group structure is induced from path multiplication. This definition requires much explication.

First, a {\it homotopy} between paths $\gamma_0$ and $\gamma_1$ is a map $H: [0,1]\times [0,1] \rightarrow E$ such that $H(s,i)=\gamma_i(s)$ for $i=0,1.$ That the homotopy fixes the boundary means that $H(0,t)=H(1,t)=*$ for all $t \in [0,1]$. We say that two paths are {\it homotopic} if there is a homotopy between them. This induces an equivalence relation on the set of paths; we call an equivalence class a {\it loop} in space. If $\alpha,\beta:([0,1],\{0,1\}) \rightarrow (E,*)$ are a pair of paths, then they can be multiplied by the rule 
$$\alpha \cdot \beta(s) = \left\{ \begin{array}{lr} \alpha(2s) & \ \  {\mbox{\rm if}}\ \  s\in[0,1/2], \\ \beta(2s-1) & \ \ {\mbox{\rm if}} \  \ s\in[1/2,1]. \end{array}\right.$$
This means that a particle traveling along the composition first travels along $\alpha$ at double speed and then along $\beta$ at double speed. 

We define a group structure on the set of homotopy classes of loops by declaring the composition to be induced by path multiplication, the identity element is the equivalence class of the constant path, and the inverse of $\alpha$ is to traverse $\alpha$ backwards (specifically $\alpha^{-1}(s)=\alpha(1-s)$). The illustrations in Fig.~\ref{fundy} indicate the geometric notions, and outline the homotopies that are needed to demonstrate that $\alpha \cdot \alpha^{-1}$ is homotopic to the constant map, and that $\cdot$ induces an associative product on equivalence classes. In this figure, the exterior of the trefoil knot is shown as the space that is interior to the torus. It is not easy to make the intuitive leap from the complement of the knot as depicted in a non-compact $3$-dimensional space to the concise picture given.

\begin{figure}[htb]
\begin{center}
\includegraphics[width=3in]{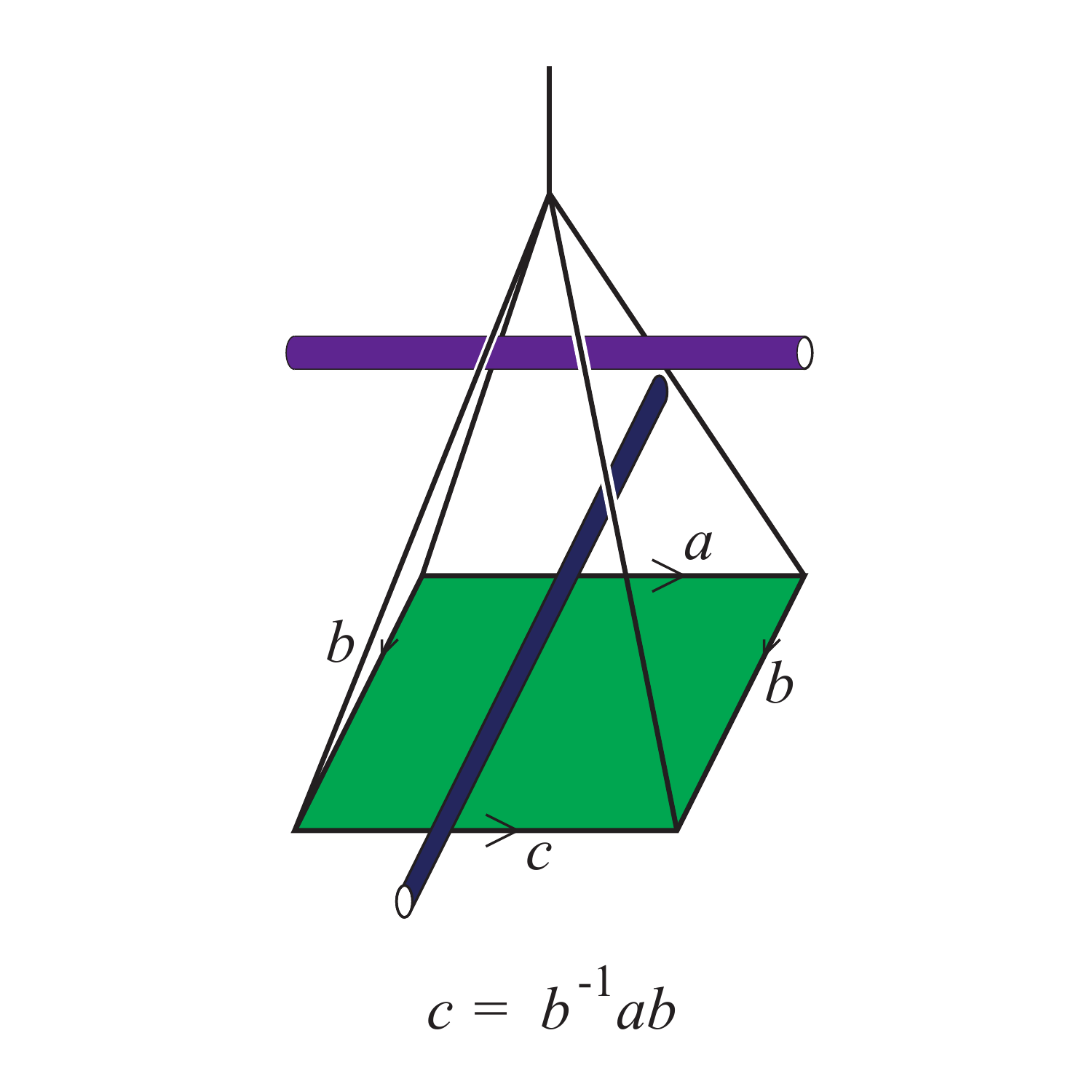}
\end{center}
\caption{The Wirtinger relations at a crossing}
\label{wirt}
\end{figure}
Following the discussion about handles in the next section, we will demonstrate that the fundamental group can be generated by the set of loops that encircle each of the arcs that appear in the diagram. We will also demonstrate that the set of relations in such a {\it knot group presentation} --- {\it i.e.} a presentation of the fundamental group of the exterior of the knot given in terms of generators and relations --- are all of the form $c=b^{-1} a b$. Figure~\ref{wirt} indicates the relationship. The triangular loop labeled $b$ on the left of the figure is homotopic to the loop labeled $b$ on the right. The square at the bottom of the figure indicates that there is a homotopy from $c$ to $b^{-1}ab$. The main portion of the homotopy is that square.

The connection between the quandle colorings given in the second section and the fundamental group is that the colorings induce representations of the fundamental group into a permutation group.  For the trefoil, the representation is into the permutations of $0,1,2$; for the figure-$8$ knot the representation is into the group of alternating permutations of $0,1,2,3$.

Figure~\ref{wirt} indicates that the relations of the form $c=b^{-1}ab$ exist at any crossing. However, it is necessary to demonstrate that such relations suffice. To this end, we decompose the knot complement into a cell complex called a handle decomposition that realizes the exterior as a thickening of a $2$-dimensional complex. 
Here are some preliminary aspects on handle theory.

\section{Handles}

For any non-negative integer $k$, let $D^k =\{ (x_1,\ldots, x_k) \in \R^k: \sum^k_{j=1} x_j^2 \le 1 \}$ denote the {\it $k$-dimensional disk.} Let $S^{k-1}=\{ (x_1,\ldots, x_k) \in \R^k: \sum^k_{j=1} x_j^2 = 1 \}$ denote the {\it $(k-1)$-dimensional sphere} which is the boundary of the disk. In case $k=0$, the  
$k$-disk is a point and its boundary is empty. Here we will concentrate on the cases in which $k\in\{0,\ldots,n\}$ and $n=1,2,3.$ 

A {\it $k$-handle} is a subset of a topological space that is homeomorphic to $D^k \times D^{n-k}$. An explicit identification between the subset and this product of disks is assumed throughout. There are several important subsets of the $k$-handle.

\begin{itemize}

\item The {\it core disk} is the subset $D^k \times \{0\}$. 

\item The {\it attaching sphere} is the subset $S^{k-1} \times \{0 \}$.

\item The {\it co-core disk} is the disk $\{0\} \times D^{n-k}.$

\item The {\it belt sphere} is the sphere $\{0\} \times S^{n-k-1}.$

\item The {\it attaching region} is $S^{k-1} \times D^{n-k}$.
\item The {\it belt region} is $D^k \times S^{n-k-1}$.

\end{itemize}

Of course, a $k$-handle is also an $(n-k)$-handle. When we consider it as such, we are {\it turning the handle decomposition upside-down}. Every handle is homeomorphic to an $n$-disk. But we are interested in the incidence relations among the handles in a given topological space (such as a surface or the complement of a knot). The incidence relations can be determined by orienting the belt sphere of a $k$-handle and computing a signed intersection between this sphere and the oriented attaching sphere of a $(k+1)$-handle. Our illustrations will attempt to make these ideas more clear. The attaching sphere and the belt sphere are often call the $A$-sphere or $B$-sphere, respectively. The basic terminology was introduced in the book \cite{RS}. 

\begin{figure}[htb]
\begin{center}
\includegraphics[width=6in]{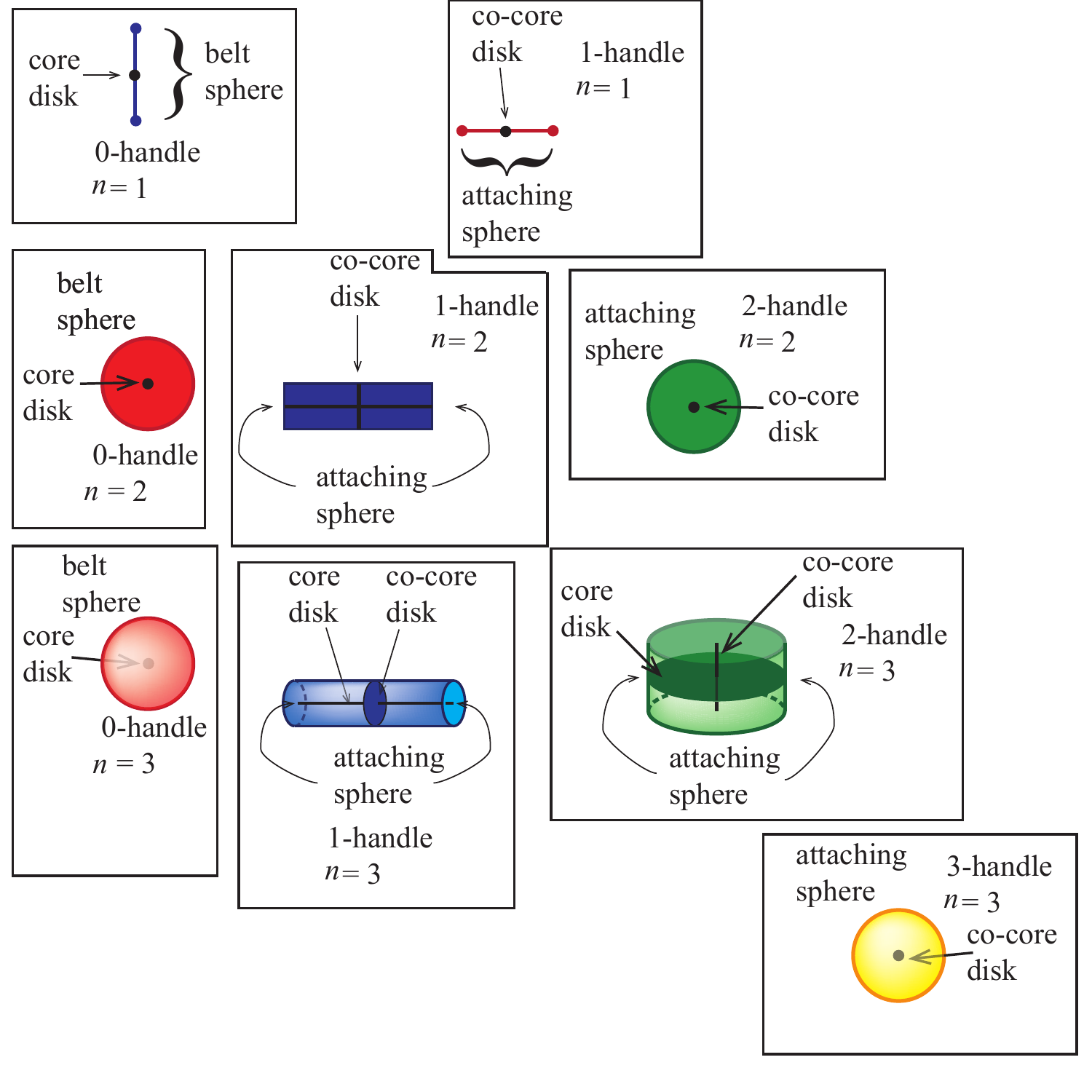}
\end{center}
\caption{Handles in dimensions $0$ through $3$}
\label{handles}
\end{figure}

Figure~\ref{handles} indicates handles in various dimensions and the important spheres and disks. Recall that a $0$-disk is a point with empty boundary.

\begin{figure}[htb]
\begin{center}
\includegraphics[width=6in]{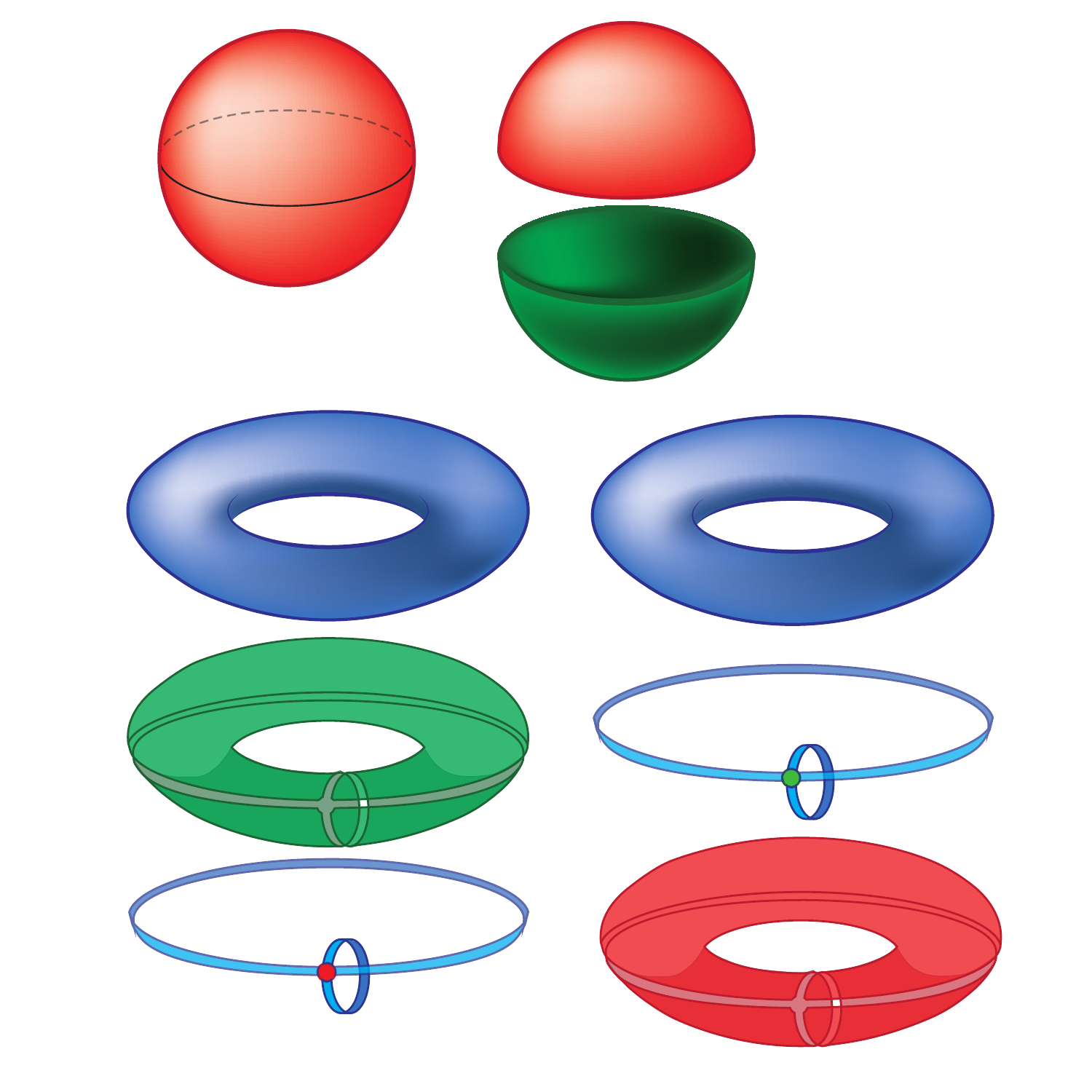}
\end{center}
\caption{Handle decompositions for the sphere and the torus}
\label{orient}
\end{figure}

In dimension $2$, compact  surfaces are classified according to their genus and the number of boundary components. The orientable surfaces that do not have boundary are in the ordered list that begins: sphere, torus, genus two surface, and that continues without end. It is instructive to give handle decompositions of some of these. Standard decompositions for the sphere and torus are given in Fig.~\ref{orient}. The sphere is decomposed as a $0$-handle and a $2$-handle. The torus is decomposed as a $0$-handle, two $1$-handles, and one $2$-handle. 
On the right of the standard handle decomposition of the torus, the same decomposition with the roles of the handles reversed is indicated. We will use this decomposition as part of the decomposition that is associated to the Alexander-Briggs presentation. In the dual decomposition the large rectangular region that wraps around most of the torus is to be considered a $0$-handle. The $1$-handles are the same, but their cores and co-cores have been switched, and the $2$-handle is the button like disk that holds the configuration together.

\begin{figure}[htb]
\begin{center}
\includegraphics[width=6in]{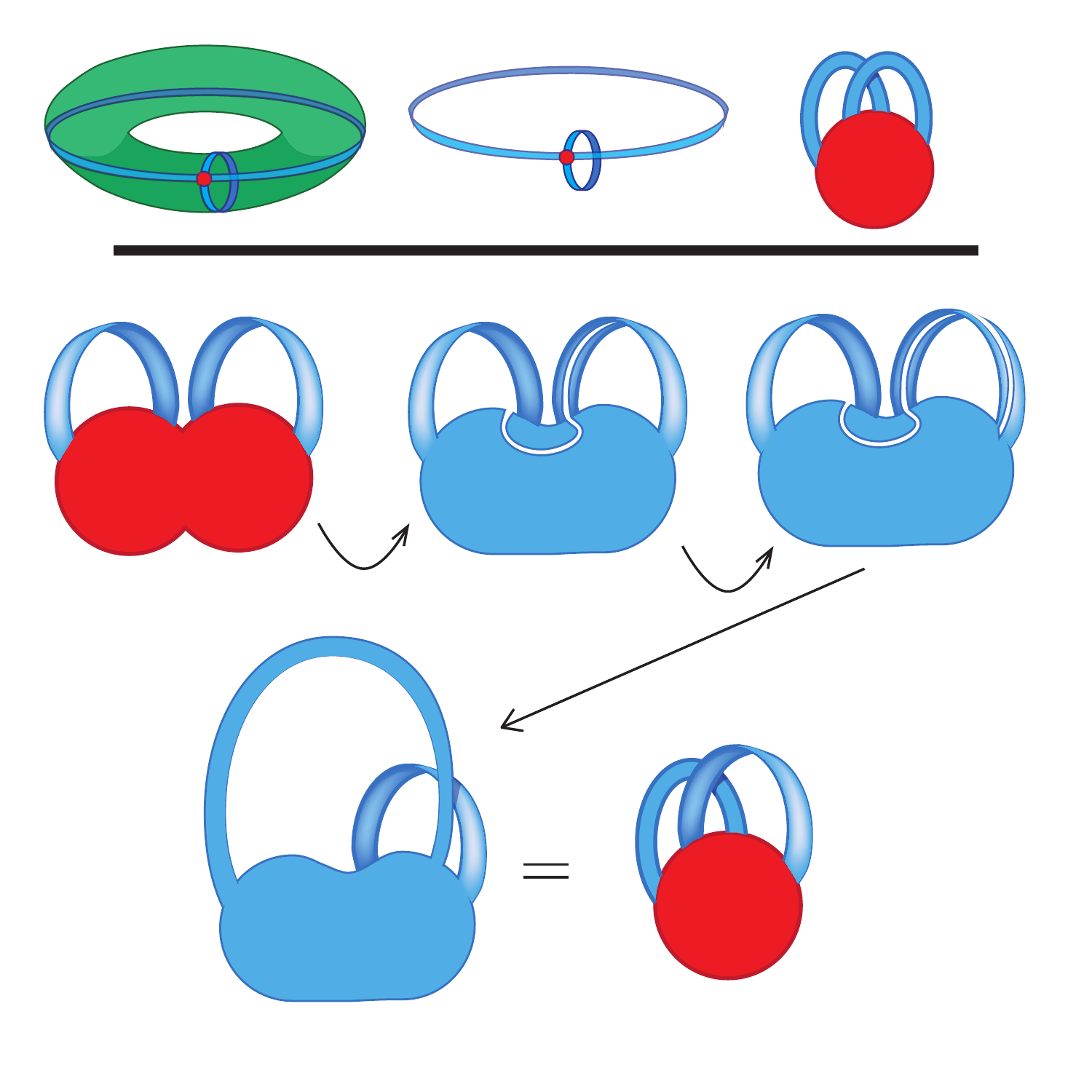}
\end{center}
\caption{The torus and two decompositions of the Klein bottle}
\label{TorusPPandKlein}
\end{figure}

Figure~\ref{TorusPPandKlein} indicates another standard picture of the handle decomposition of the torus. This time the $2$-handle has been removed, and the meridional and longitudinal $1$-handles appear as ``basket" handles attached to the disk. In the lower figure, the Klein bottle is indicated with a $2$-handle missing, and it has been constructed as a schematic of the connected sum of two projective planes. The projective planes having been obtained as a single $0$-handle, a $1$-handle attached with a twist, and the $2$-handle to be attached on the outside. In the next two illustrations a topological deformation occurs in which the M\"obius band on the left ``slides over" the M\"obius band on the right. The sliding is achieved via the method of excavating a thin disk from the surface. At the bottom left of the illustration the resulting $1$-handle is bent upwards and rounded out. The resulting decomposition is an untwisted band nestled with a twisted band as illustrated on the bottom right.

\begin{figure}[htb]
\begin{center}
\includegraphics[width=6in]{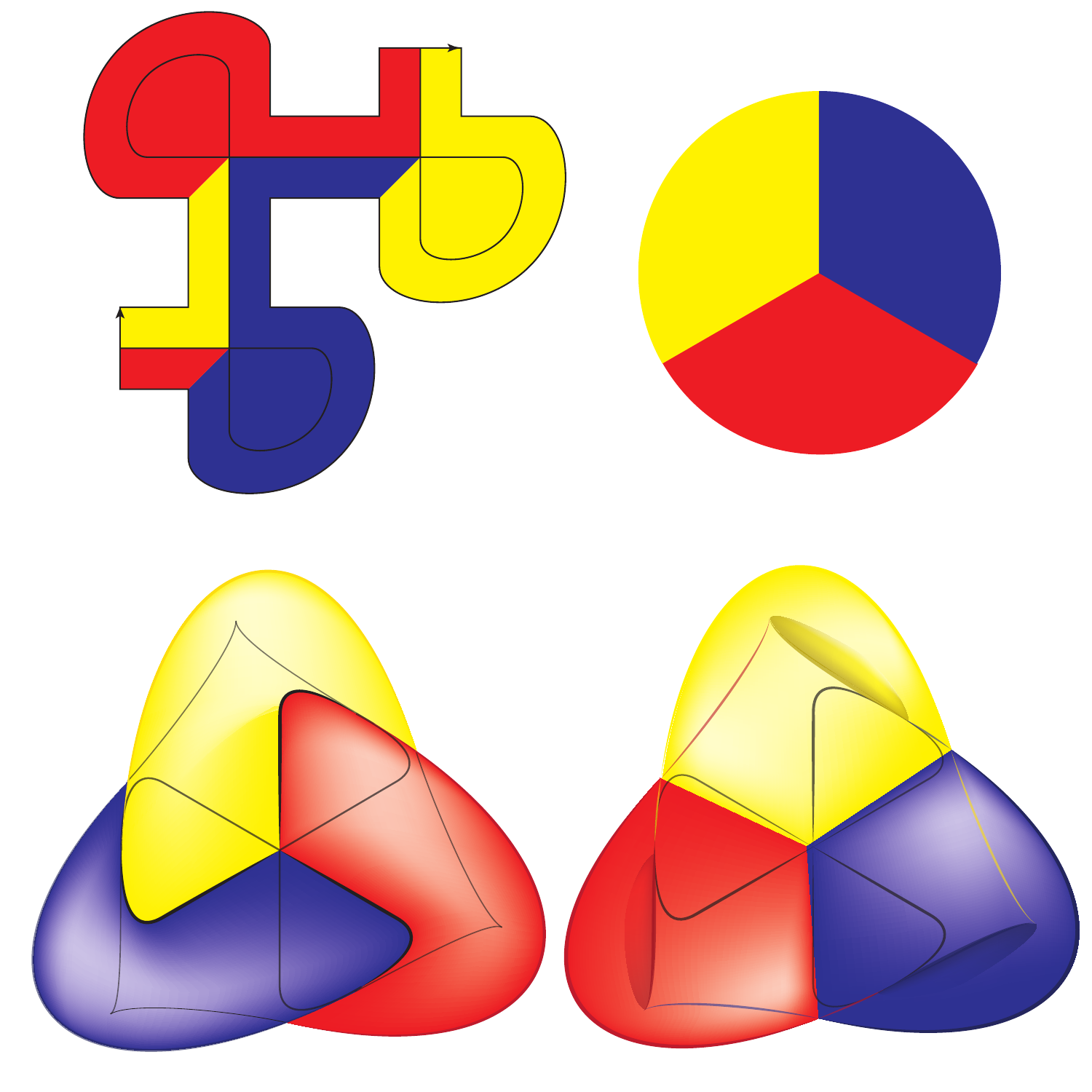}
\end{center}
\caption{Boy's surface, front, back, and as a handle decomposition}
\label{newboy}
\end{figure}

\begin{figure}[htb]
\begin{center}
\includegraphics[width=2in]{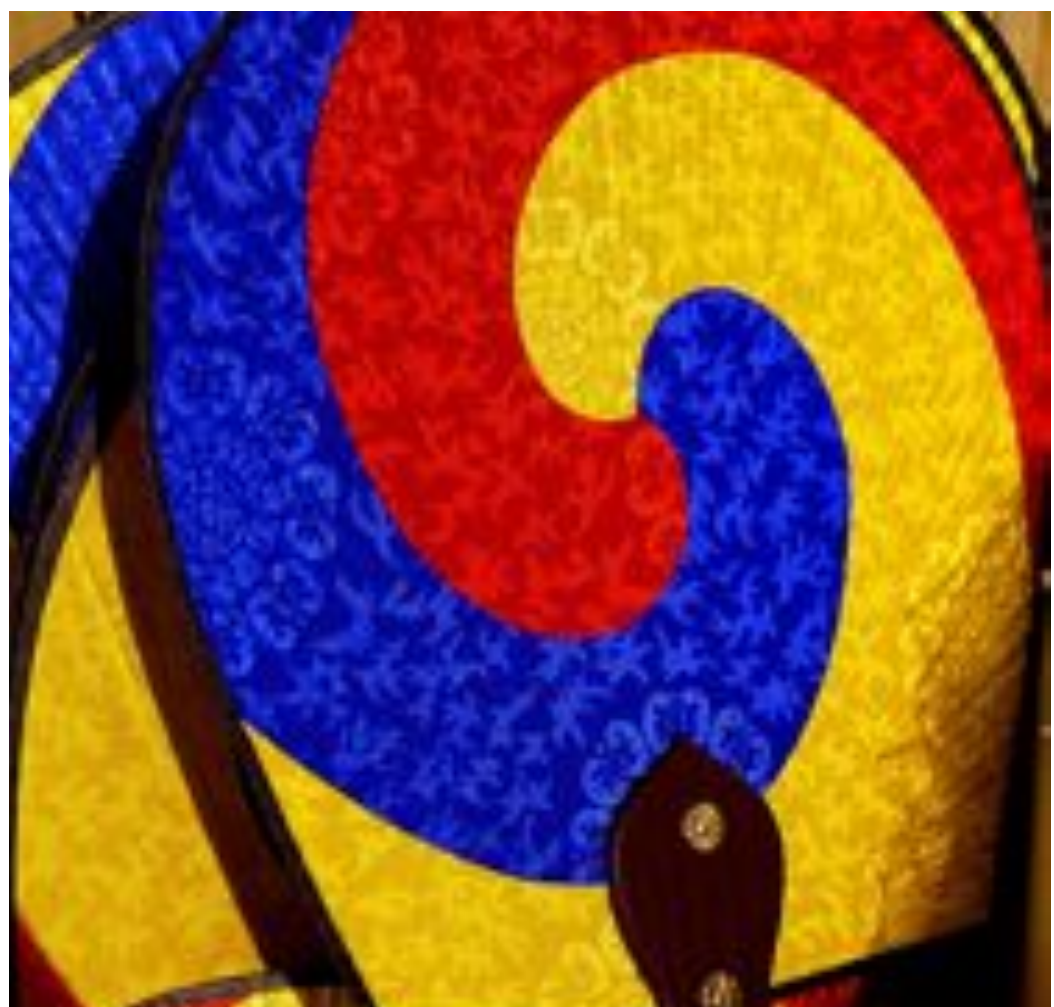}
\end{center}
\caption{Traditional Korean fans}
\label{fan}
\end{figure}

Boy's surface (Fig.~\ref{newboy}) is an immersion of the projective plane that is immersed in $3$-dimensional space. It was discovered by David Hilbert's student Werner Boy much to the surprise of Hilbert. It enjoys many remarkable properties including a $3$-fold symmetry. The symmetry here is indicated  by means of primary colors: red, blue, and yellow. These were chosen to emulate the colors and symmetries expressed in some traditional Korean fans (Fig.~\ref{fan}) in homage to my hosts.

\section{Handles in dimension $n=3$}

Handle techniques are often used in the study of $3$ and $4$ dimensional manifolds. In $3$-dimensions, the handle decomposition coincides with what is also called the Heegaard decomposition of the $3$-dimensional manifold. 
As an example, the complement of a knot is a $3$-dimensional manifold with boundary. We will use the knot diagram to construct a handle decomposition of the knot complement. In this decomposition, the arcs correspond to $1$-handles and the crossings (under arcs) correspond to $2$-handles. 

\begin{figure}[htb]
\begin{center}
\includegraphics[width=3in]{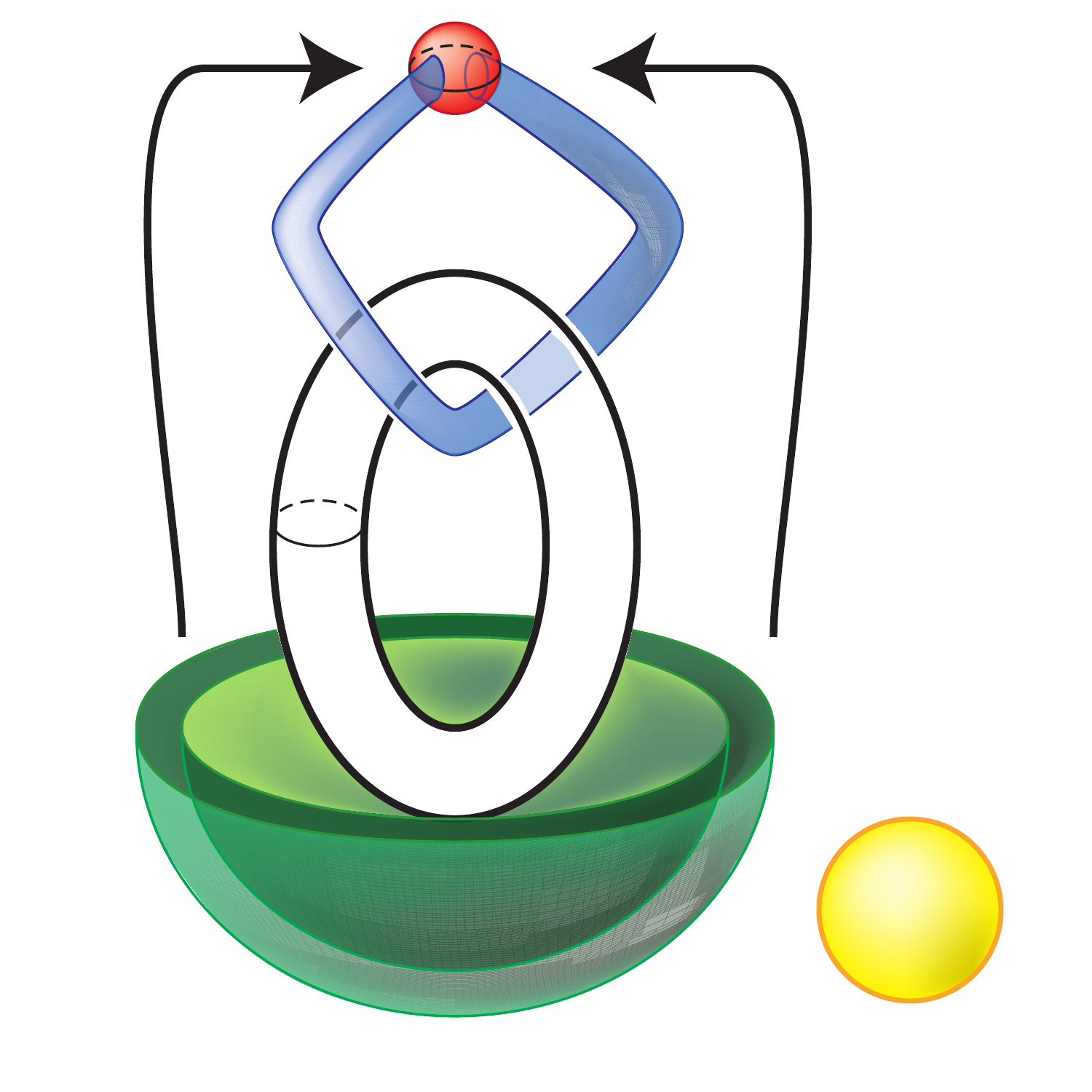}
\end{center}
\caption{A handle decomposition of the unknot}
\label{handledecomunk}
\end{figure}

In order to help develop this description in general, we begin with a decomposition of the unknot (Fig.~\ref{handledecomunk}). A $0$-handle appears at the top of the diagram. A $1$-handle is attached in such a way that its core disk corresponds to the disk that is bounded by the circle in space. The $1$-handle, then corresponds to a maximal point of the curve as it appears in space. A $2$-handle is attached at the minimal point of the curve. The attaching sphere for the $2$-handle does not intersect the core disk of the $1$-handle. A $3$-handle can be attached on the outside to complete the construction of the complement.

\begin{figure}[htb]
\begin{center}
\includegraphics[width=6in]{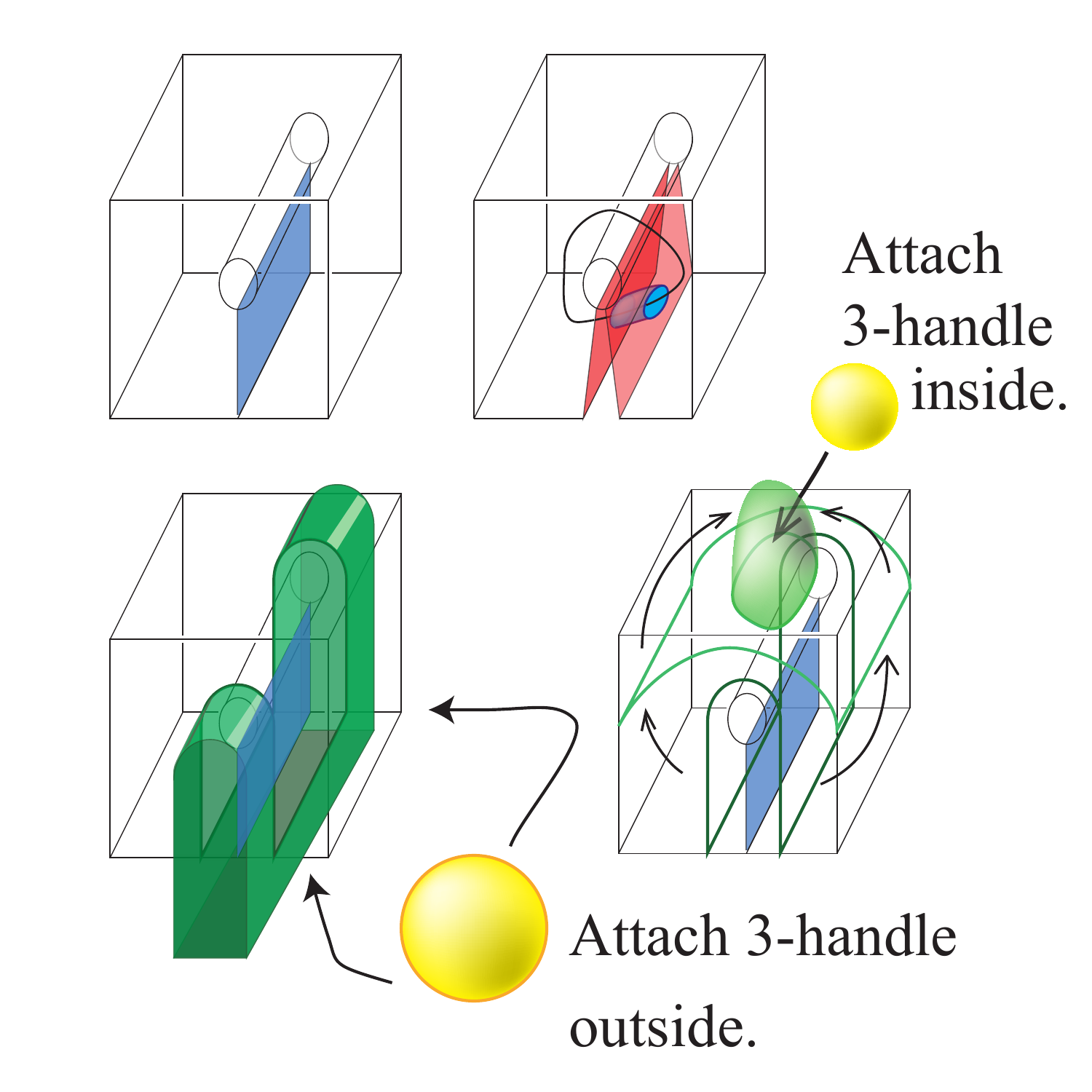}
\end{center}
\caption{Re-interpreting this handle decomposition of the unknot}
\label{keyhole}
\end{figure}

To better understand this decomposition (and those which follow) the picture will be rearranged. Consider first the complement of a neighborhood of the maximal point of the circle. It is homeomorphic to a cube with a hole drilled through it as in the left of Fig.~\ref{keyhole}. Observe that the
the cube from which a hole has been drilled is homeomorphic to a $3$-ball to which a $1$-handle has been attached. Figure~\ref{keyhole} indicates the correspondences. Next observe that the attaching sphere for the $2$-handle (which corresponds to the minimum point in the diagram) does not intersect the belt sphere of the $1$-handle. The core disk of the $2$-handle envelops the two ends of the excavated hole. The $3$-handle consists of the space in which the observer sits. 

Alternatively, the attaching sphere for the $2$-handle can be slid to the top of the cube in the figure, and the $3$-handle is nestled inside the resulting purse. 

\section{Handles in $3$-dimensions}

A knot diagram easily lends itself to the construction of a handle decomposition of the complement. The arcs of the diagram correspond to holes that have been drilled from the space above the diagram as indicated in Fig.~\ref{keyhole}. Such ``holes removed" can be reinterpreted as $1$-handles that are attached to a $0$-handle that is envisioned as lying above the plane in which the diagram is drawn. The crossings --- representing arcs that go under the plane of the diagram --- define $2$-handles. 

\begin{figure}[htb]
\begin{center}
\includegraphics[width=6in]{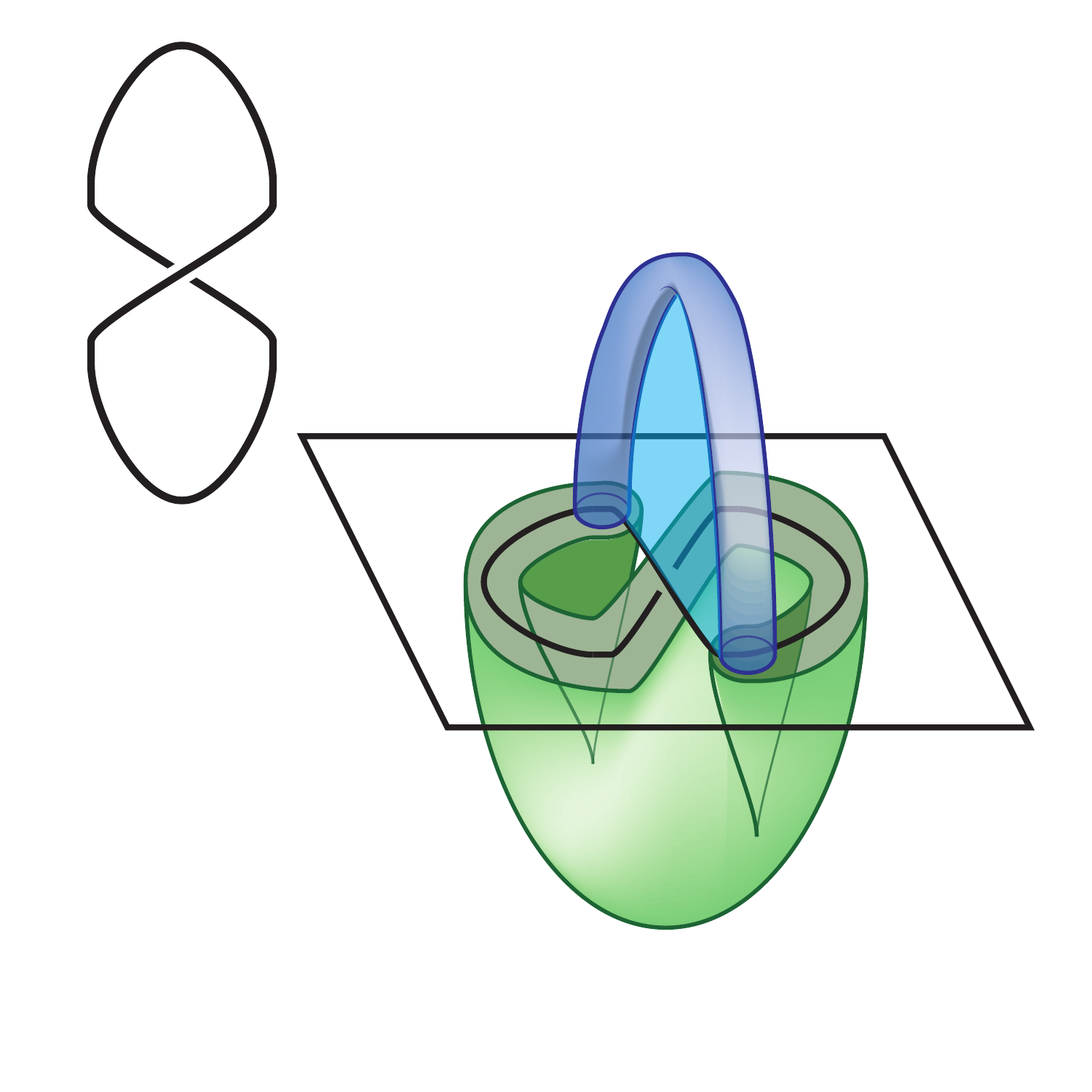}
\end{center}
\caption{Another handle decomposition of the unknot}
\label{other8}
\end{figure}

The most simple case that can be depicted and that includes a crossing is that of the unknot with one crossing. The handle decomposition is represented in Fig.~\ref{other8}. This figure indicates that the attaching sphere for the $2$ handle (which is an $S$-shape in the figure) intersects the belt sphere of the $1$-handle in a sequence $1, -1,-1,1$. Here there is one $1$-handle (thus labeled ``1") and is oriented from the lower right towards the upper left of the planar picture, and the  $S$-curve is oriented clockwise. Clearly, the intersections between the $1$ and $2$ handles can be removed since there are a pair of disks co-bounded by segments of each of these handles.

\begin{figure}[htb]
\begin{center}
\includegraphics[width=6in]{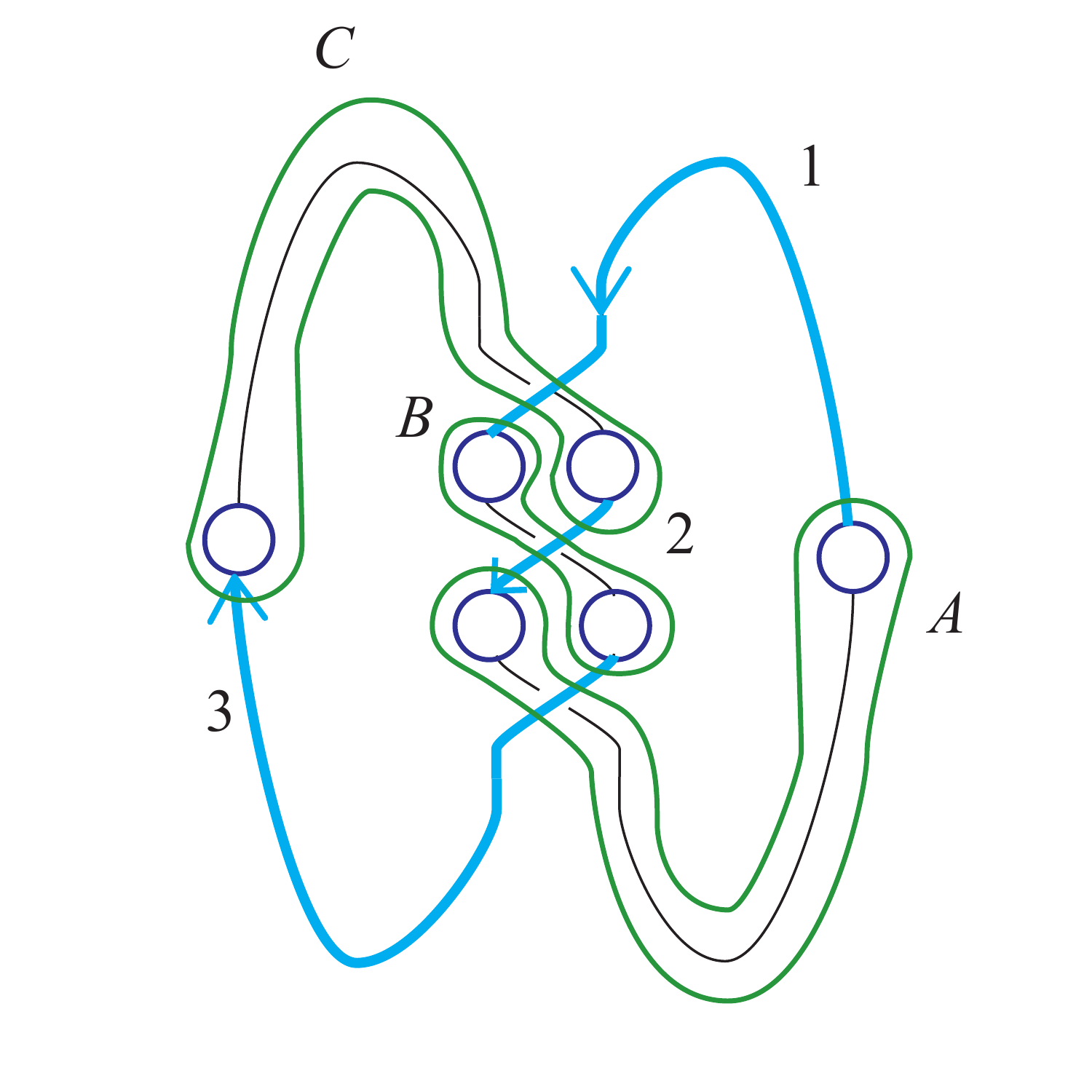}
\end{center}
\caption{A handle decomposition of the trefoil knot}
\label{trehandle}
\end{figure}

The intersections between $1$-handles and $2$-handles can be indicated completely within a planar diagram built from the knot diagram. Figure~\ref{trehandle} demonstrates one such diagram associated to the trefoil. Over arcs 
correspond to the segments of the co-cores of the $1$-handles that intersect the plane of the diagram. The attaching spheres for the $2$-handles encircle the pair of segments of the under crossing.

In the diagram of Fig.~\ref{trehandle}, we can read the intersection of the attaching sphere of the $2$-handles with the co-cores of the $1$-handles. In this diagram, the attaching spheres for the $2$-handles are oriented counterclockwise. So the intersection sequences are (up to cyclic permutation) as follows:
$$A: 1,3^{-1},2^{-1},3;$$
$$B: 3,2^{-1},1^{-1},2;$$
$$C: 2,1^{-1},3^{-1},1.$$

Any path in space is homotopic to a path that misses the $3$-handle. Thus any loop in space can be moved into the union of the $0$, $1$, and $2$-handles. We realize, therefore that the intersection sequences $A$,$B$, and $C$ define a presentation for the fundamental group of the trefoil.
Specifically, the fundamental group is given:
$$\langle x,y,z : xz^{-1}y^{-1}z, \ \ zy^{-1}x^{-1}y, \ \ yx^{-1}z^{-1}x \rangle.$$

\begin{figure}[htb]
\begin{center}
\includegraphics[width=6in]{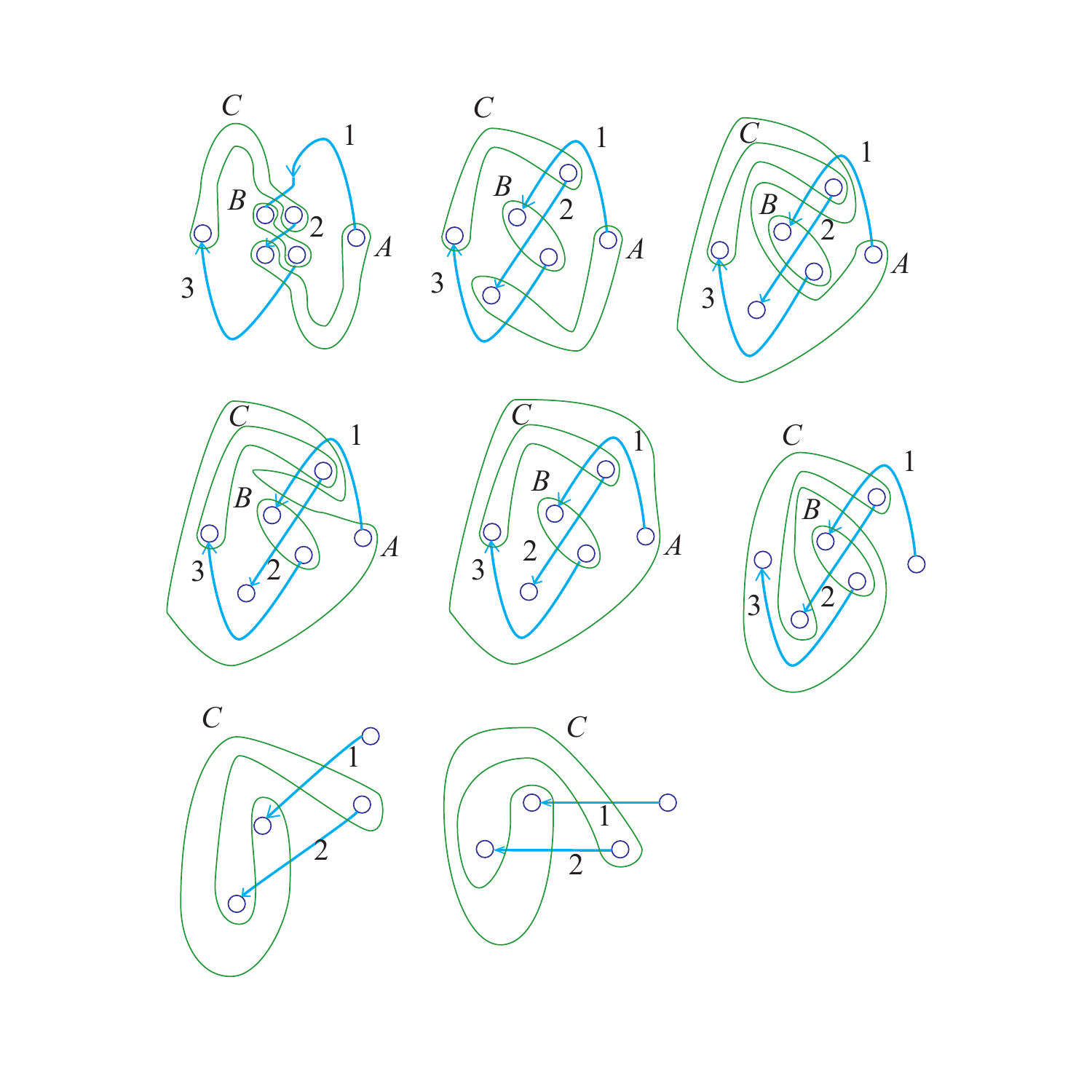}
\end{center}
\caption{Handle sliding and cancelation of the trefoil complement}
\label{trehandleslides}
\end{figure}

Just as in the case of surfaces, handles in the decompositions of knot exteriors (or any $3$-dimensional manifold) can be slid over each other. Moreover, it is possible to cancel $i$-handles against $(i+1)$-handles provided they intersect geometrically once. A sequence of handle slides and cancelations is depicted in Fig.~\ref{trehandleslides}. Let the diagrams in the figure be labeled by row/column indices. The transition from $(1,1)$ to $(1,2)$
is simply a topological isotopy. From $(1,2)$ to $(1,3)$  handle $A$ slides over handle $C$. From $(1,3)$ to $(2,1)$ handle $A$ slides over handle $B$. From $(2,1)$ to $(2,2)$ handle $A$ is homotoped to have no intersections with the $1$-handles. At this point handle $A$ can cancel with the $3$-handle on the outside. The move from $(2,2)$ to $(2,3)$ is this cancelation together with handle $C$ sliding over $B$. From $(2,3)$ to  $(3,1)$ handle $3$ cancels with $B$. The move from $(3,1)$ to $(3,2)$ is a homotopy. 

The presentation for the fundamental group reduces to 
$$\langle x ,y : yx^{-1}y^{-1}x^{-1}yx \rangle = 
\langle x ,y : xyx = yxy \rangle.$$

\begin{figure}[htb]
\begin{center}
\includegraphics[width=6in]{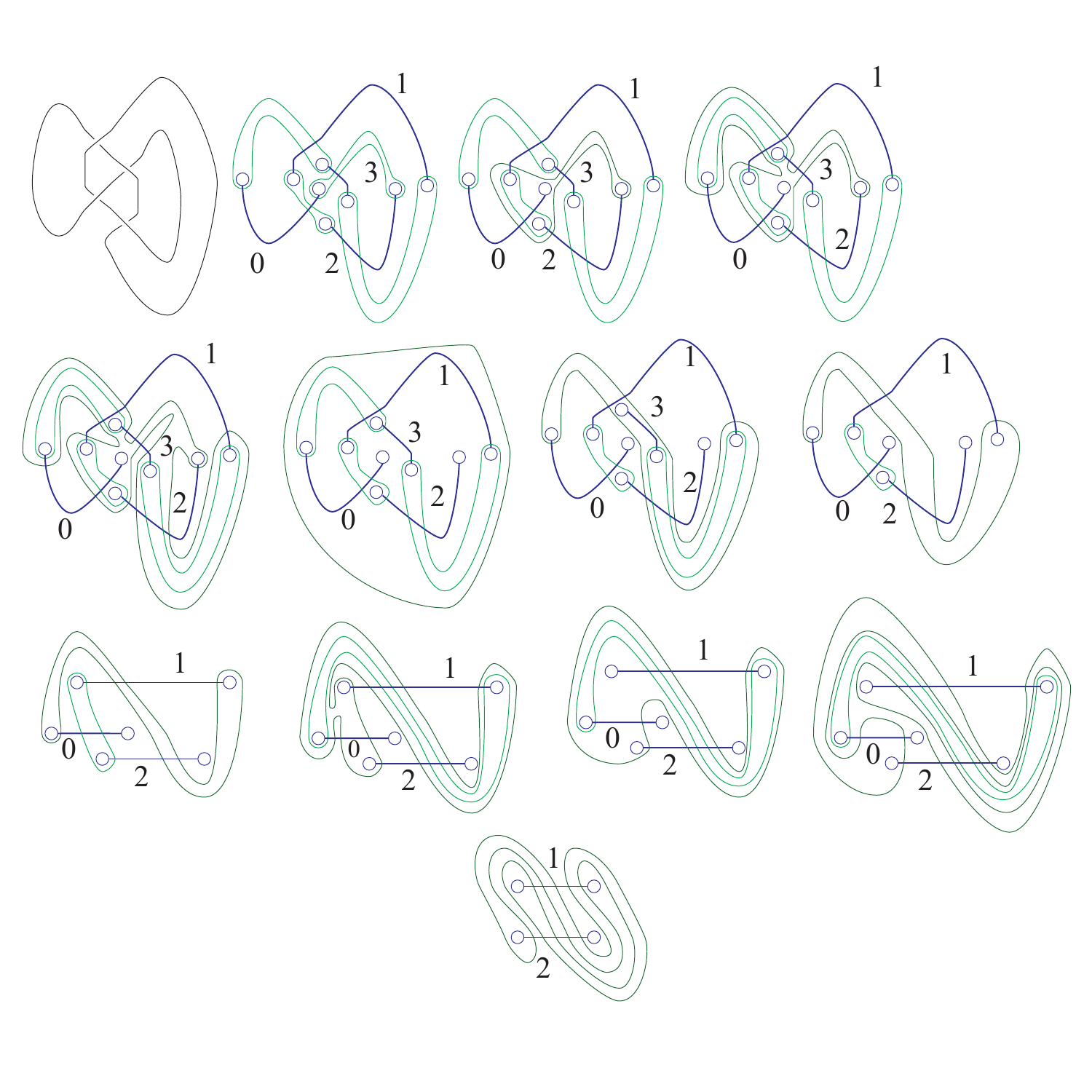}
\end{center}
\caption{Handle sliding on the figure-$8$ knot complement}
\label{atewirt}
\end{figure}

Figure~\ref{atewirt} indicates a similar sequence of handle slides that moves the handle decomposition of the figure-8 knot to a decomposition with one $2$-handle and two $1$-handles. The last illustration in this sequence is related to the so-called $2$-bridge presentation of the knot.  From this representation one can determine that the $2$-fold branched cover is the $3$-dimensional manifold which is known as the lens space $L(5,2)$: this manifold has a handle decomposition with one $0$-handle, one $1$-handle and one $2$-handle that is wrapped around the boundary torus twice longitudinally and five times meridonally. 

Some general properties can be proven following a close examination of these examples. First of all, it is always the case that one of the $2$-handles can be cancelled against the $3$-handle on the outside. This can be seen by mimicking the handle slides in these examples. Any particular $2$-handle's attaching sphere can be slid to envelope the remaining attaching spheres. Consequently, the knot complement is, in essence, a $2$-dimensional complex that can be formed by attaching as many $1$-handles as over-arcs, and one fewer $2$-handle. 

Often the handle decomposition can be simplified. In particular, whenever the attaching sphere for a $2$-handle intersects the belt sphere for a $1$-handle exactly once and this belt sphere has no other intersections, then this pair of handles can be removed. The removal corresponds to removing a generator of the fundamental group and replacing all occurrences of that generator with a sequence that is obtained from intersection sequence of the attaching sphere of the $2$-handle.

All of the aspects of group presentations for classical knots and the relationship between this presentation (the Wirtinger presentation) and the handle decomposition induced from the under-arcs can be found in various spots in the literature.

\begin{center}\rule{3in}{0.005in}\end{center}

There are many other methods to obtain a presentation for the fundamental group of a knot complement, and even though many important knot invariants can be obtained from the fundamental group (and the closely related fundamental quandle), other invariants such as the Jones polynomial \cite{Jones} and the HOMPLYPT (or FLY THOMP) \cite{homfly} polynomial are obtained diagrammatically, but their relationship to other traditional topological invariants is difficult to discern. The is not enough time or space within this article to spend time on these. However, there is an alternative method for obtaining a presentation of the fundamental group that will close the article.

\section{The Alexander-Briggs presentation}
The {\it Alexander-Briggs} presentation of the fundamental group is obtained from a handle decomposition of the knot complement in which most (all but two) of the  $1$-handles are associated to the crossings, and  most of the $2$-handles (all but one) are associated to the regions. In future work with Masahico Saito, Dan Silver, and Susan Williams, I  will develop this presentation in the context of virtual knots and their Alexander invariants. 

\begin{figure}[htb]
\begin{center}
\includegraphics[width=4in]{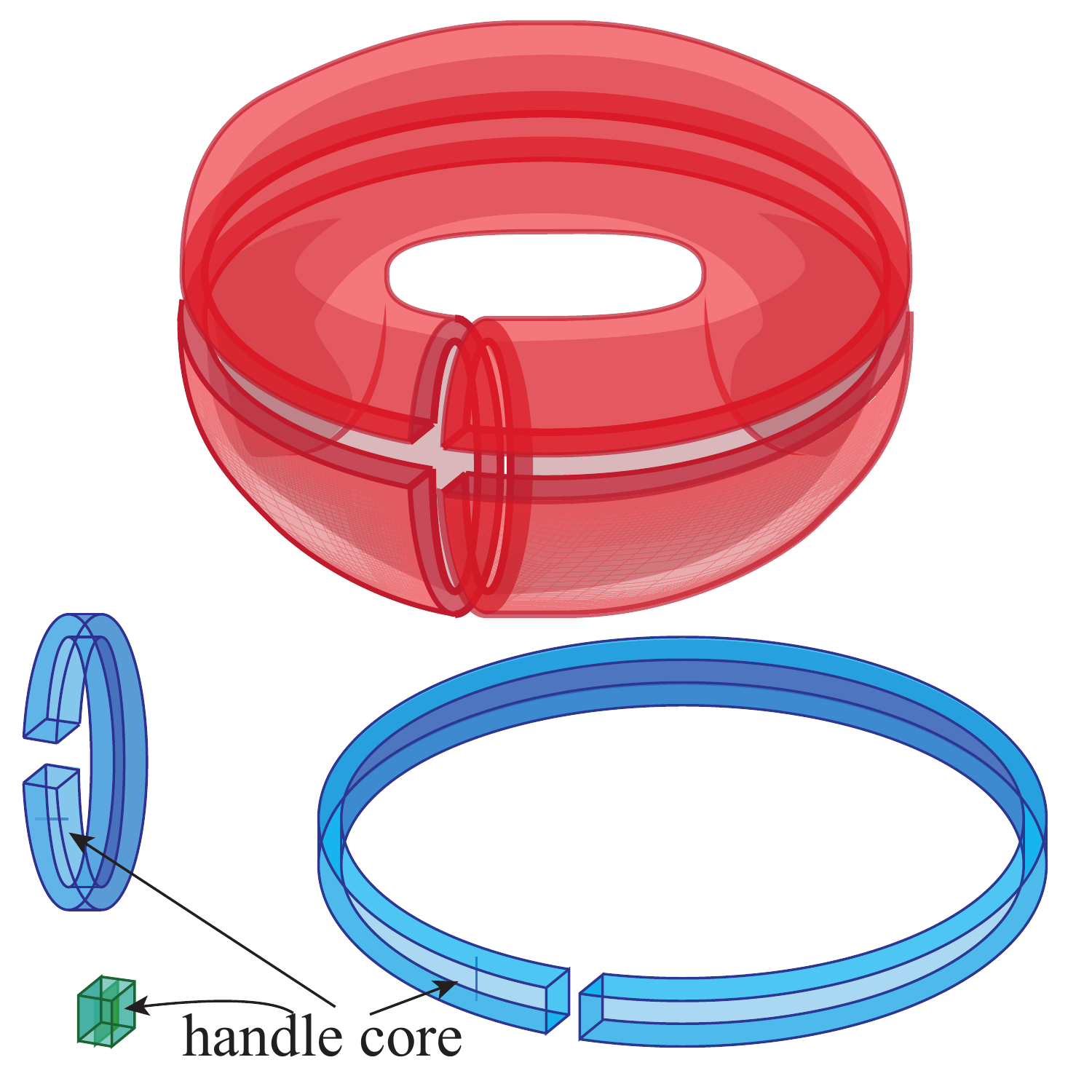}
\end{center}
\caption{The decomposition of the knot complement in the neighborhood of the boundary}
\label{ABtorus}
\end{figure}

In the initial stage of this decomposition, the complement is built out from the boundary torus. The ``upside-down" decomposition of the torus that is depicted on the lower right side of Fig.~\ref{orient} is thickened to a decomposition of the torus times an interval. The $0$-handle envelopes the torus. The $1$-handles are spatial regular neighborhoods of the letter {\sf C}, with the co-core disks being planar neighborhoods. These $1$-handles correspond to the longitude and meridian of the knot, but be aware that since the cores are short arcs, that which is apparently a longitude is a segment of the co-core disk of the meridian and {\it vice versa}. Figure~\ref{ABtorus} contains the details.

\begin{figure}[htb]
\begin{center}
\includegraphics[width=2in]{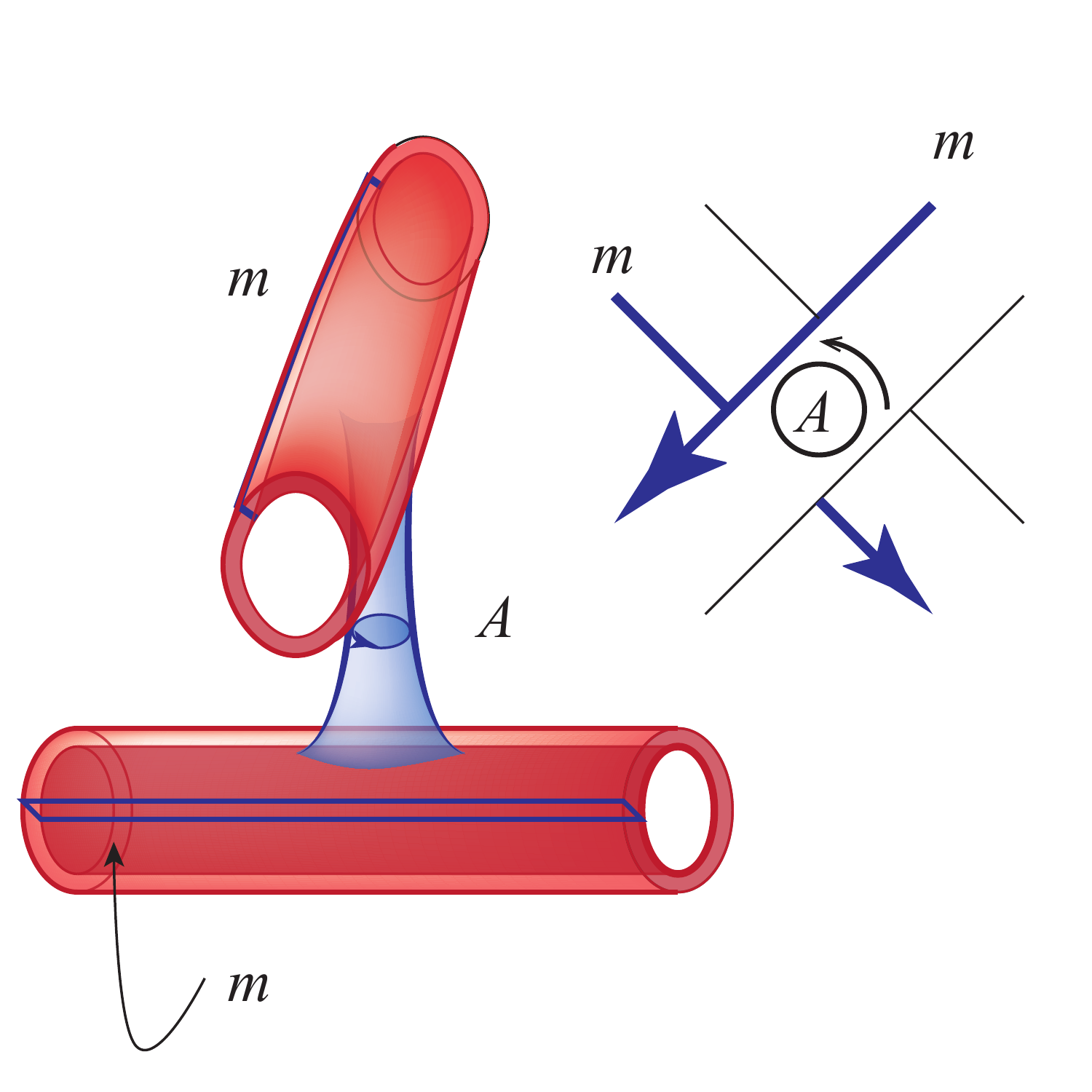}
\end{center}
\caption{At each crossing a $1$-handle is attached to the $0$-handle that envelopes the boundary}
\label{ABcrossing}
\end{figure}

At each crossing a $1$-handle is attached. It is a pillar erected between neighborhoods of arcs at the crossing. The core disk runs from the lower arc to the upper arc, and the belt sphere is oriented in a counter clockwise fashion. The attaching sphere lies on the $0$-handle that envelopes the boundary torus. Figure~\ref{ABcrossing} indicates the details.

\begin{figure}[htb]
\begin{center}
\includegraphics[width=2in]{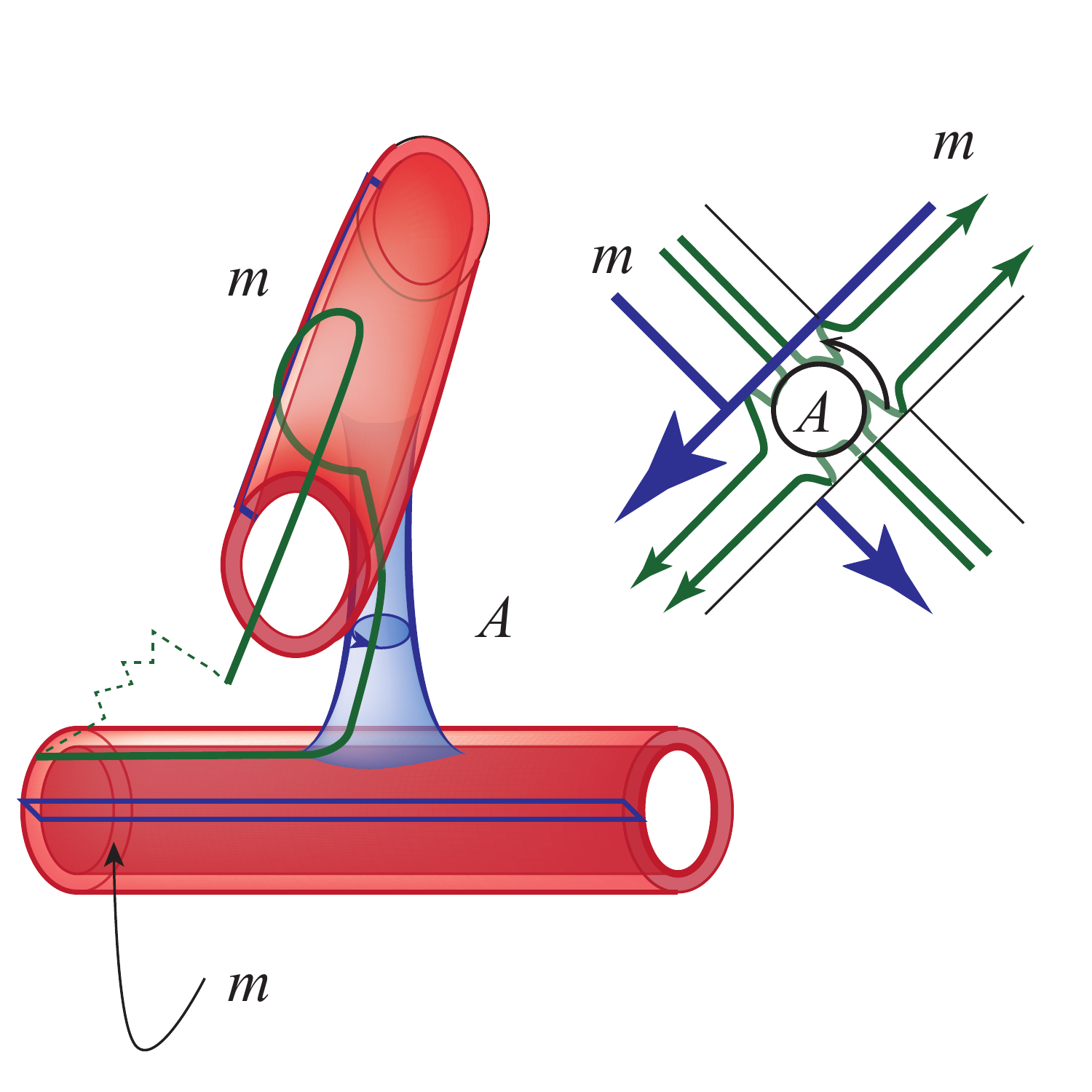}
\end{center}
\caption{The attaching spheres of the four $2$-handles that are incident to the pillar are shown in the associated schematic}
\label{ABtwohandles}
\end{figure}

There is a $2$-handle that lies in a neighborhood of the boundary. In Fig.~\ref{ABtorus} this is represented as a small cube. Its attaching sphere intersects the belt sphere of the longitude and meridian in the commutator $MLM^{-1}L^{-1}$ up to orientation and cyclic permutation. 

The remaining $2$-handles are represented by the regions in the diagram. See Fig.~\ref{fig8} for the definition of regions. The attaching sphere for the regions intersects the pillars at the crossings. Figure~\ref{ABtwohandles} indicates the intersection of one such attaching sphere in a neighborhood of the crossing. A schematic diagram on the right side of the illustration indicates how this attaching sphere  and those of the adjacent regions can be drawn at a knot diagram. In the schematic, the co-core disk of the meridian is illustrated as thick segments. Each region including the unbounded region of a planar diagram serves as the core disk of a $2$-handle.

To complete the construction of the complement a pair of $3$-handles are attached to the complex --- one $3$-handle above the plane of the diagram and one $3$-handle below the plane of the diagram.

\begin{figure}[htb]
\begin{center}
\includegraphics[width=3in]{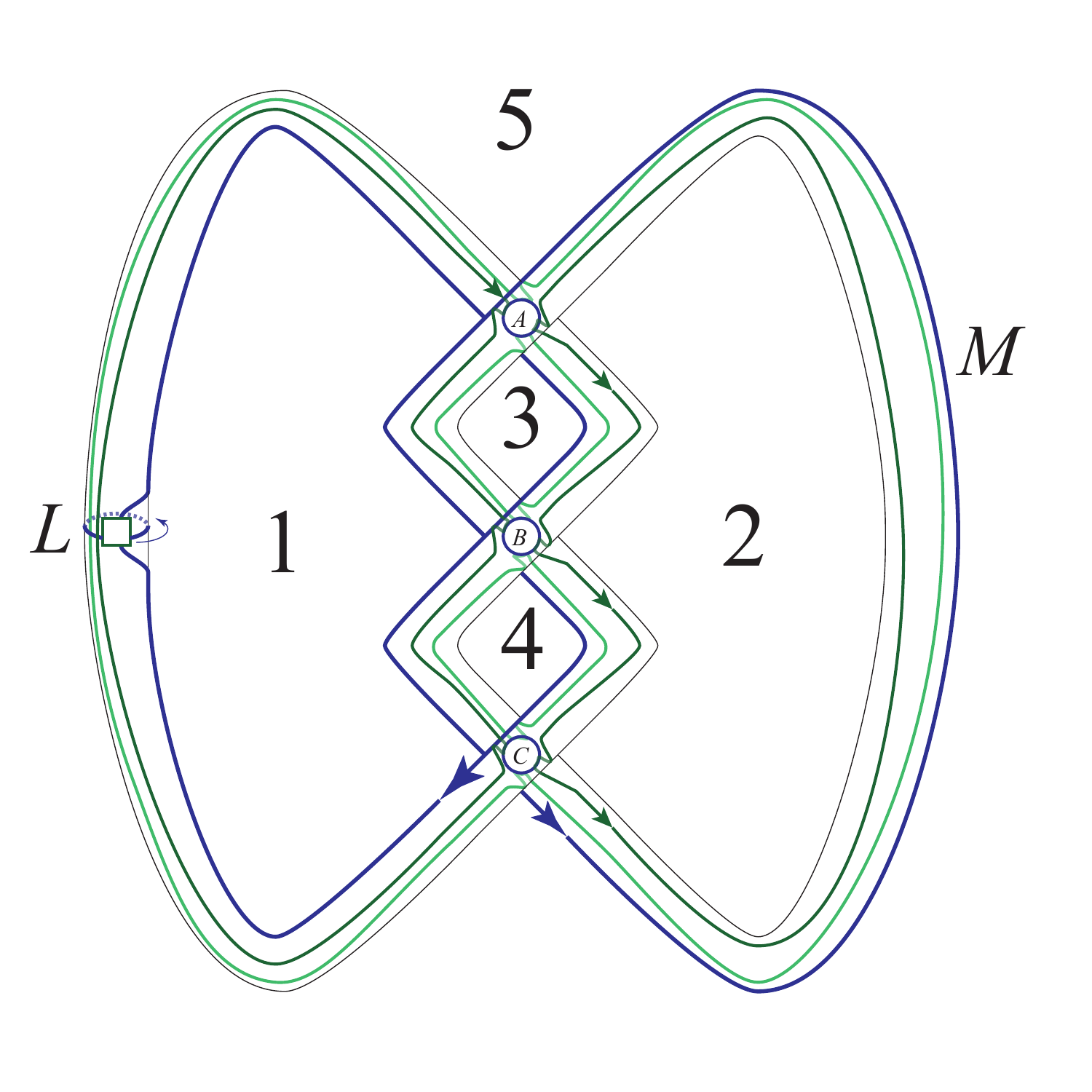}
\end{center}
\caption{The decomposition of the trefoil complement}
\label{ABtref}
\end{figure}

Figure~\ref{ABtref} indicate the handle decomposition of the trefoil knot. From this decomposition we can deduce a presentation for the fundamental group with generators 
$m$, $L$, $A$, $B$, and $C$. There are six relations that are read from the attaching spheres of the $2$-handles. These are: \begin{enumerate}
\item[0.] $L^{-1} M^{-1} L M =1.$
\item[1.] $LAMBMCM=1.$
\item[2.] $CBA=1.$
\item[3.] $ABM=1$.
\item[4.] $BCM=1.$
\item[5.] $LAMC=1$.
\end{enumerate}

From the relations (5) and (1), we obtain
$C^{-1}BMCM=1$. From (3) and (4) we obtain
$AM^{-1}C^{-1}M=1.$ Thus $C=MAM^{-1}$, and $B=A^{-1}M^{-1}$.  Substituting, 
$$C^{-1}BMCM= [ MA^{-1}M^{-1} ] [ A^{-1}M^{-1} ] M [  MAM^{-1}  ]M =1;$$
 $$MA^{-1}M^{-1} A^{-1} MA =1;$$
 $$A^{-1}M^{-1}A^{-1}=M^{-1}A^{-1}M^{-1}.$$
Consequently, the group is presented as 
$$\langle A,M: AMA=MAM \rangle.$$ 
Since the longitude can be written as $L=C^{-1}M^{-1}A^{-1}$, further substitutions give that the relation $LML^{-1}M^{-1}=1$ is a consequence of the  defining relation $AMA=MAM$. 

We close this discussion with the remark that the dot notation in \cite{AB} can be chosen to coincide with the incidence of the meridional co-core and the regions with some crossing conventions. Thus there is a calculus to go from this presentation of the fundamental group to a presentation of the Alexander module. This topic and more will be the subject of a subsequent paper.

\bibliographystyle{mdpi}
\makeatletter
\renewcommand\@biblabel[1]{#1. }
\makeatother

\end{document}